\documentclass[11pt]{article}

\def\mathbb{\Bbb}
\usepackage{amsfonts,latexsym,amstext}
\usepackage{amsmath}
\usepackage[active]{srcltx}

\usepackage{amssymb}
\usepackage{amsmath,amscd,amssymb}

\newtheorem{theorem}{Theorem}[section]
\newtheorem{lemma}[theorem]{Lemma}
\newtheorem{proposition}[theorem]{Proposition}
\newtheorem{definition}{Definition}[section]

\newtheorem{hypothesis}[theorem]{Hypothesis}

\newtheorem{remark}[theorem]{Remark}

\numberwithin{equation}{section}

\def\qed{{\hfill\hbox{\enspace${ \square}$}} \smallskip}
\def\sqr#1#2{{\vcenter{\vbox{\hrule height .#2pt \hbox{\vrule
 width .#2pt height#1pt \kern#1pt \vrule
width .#2pt} \hrule height .#2pt}}}}
\def\square{\mathchoice\sqr54\sqr54\sqr{4.1}3\sqr{3.5}3}

\def\ds{\begin{displaystyle}}
\def\eds{\end{displaystyle}}
\def\dis{\displaystyle }
\def\<{\langle }
\def\>{\rangle }

\def\R{\mathbb R}
\def\N{\mathbb N}

\def\E{\mathbb E}
\def\P{\mathbb P}
\def\Q{\mathbb Q}

\def\calf{{\cal F}}
\def\calg{{\cal G}}

\def\calm{{\cal M}}
\def\caln{{\cal N}}

\def\call{{\cal L}}
\def\cals{{\cal S}}

\newcommand{\abs}[1]{\left\vert#1\right\vert}

\title{HJB equations in infinite dimension with locally Lipschitz Hamiltonian
and unbounded terminal condition}

\begin{document}

\title{
}
\date{}
 \author{
 Federica Masiero\\
 Dipartimento di Matematica e Applicazioni, Universit\`a di Milano Bicocca\\
 via Cozzi 55, 20125 Milano, Italy\\
 e-mail: federica.masiero@unimib.it,\\
\\
Adrien Richou\\
Univ. Bordeaux, IMB, UMR 5251, F-33400 Talence, France.\\
CNRS, IMB, UMR 5251, F-33400 Talence, France.\\
INRIA, \'Equipe ALEA, F-33400 Talence, France.\\ 
e-mail: adrien.richou@math.u-bordeaux1.fr}

\title{HJB equations in infinite dimensions with locally Lipschitz Hamiltonian
and unbounded terminal condition}

\maketitle  
\begin{abstract}
We study Hamilton Jacobi Bellman equations in an infinite dimensional Hilbert space,
with Lipschitz coefficients, where the Hamiltonian has superquadratic growth with
respect to the derivative of the value function, and the final condition is not bounded.
This allows to study stochastic optimal control problems for suitable
controlled state equations with unbounded control processes. The results are applied to a
controlled wave equation.
\end{abstract}

\section{Introduction}

In this paper we study semilinear Kolmogorov equations
in an infinite dimensional Hilbert space $H$ of the following form:
\begin{equation}
\left\{
\begin{array}
[c]{l}%
\frac{\partial v}{\partial t}(t,x)=-\mathcal{L}v\left(  t,x\right)
+\psi\left( t,x,v(t,x), \nabla v\left(  t,x\right)B  \right)  ,\text{ \ \ \ \ }t\in\left[  0,T\right],
\text{ }x\in H,\\
v(T,x)=\phi\left(  x\right).
\end{array}
\right.  \label{Kolmo intro}%
\end{equation}
$\mathcal L $ is the generator of the transition
semigroup $P_t$ related to the following perturbed Ornstein-Uhlenbeck process in $H$
\begin{equation}
\left\{
\begin{array}
[c]{l}%
dX_\tau^{t,x}  =AX_\tau^{t,x} d\tau+B G(X_\tau^{t,x}) d\tau+BdW_\tau
,\text{ \ \ \ }\tau\in\left[  t,T\right], \\
X_t^{t,x} =x,
\end{array}
\right. \label{ornstein-pert intro}%
\end{equation}
where $W$ is a cylindrical Wiener process with values in another real and separable
Hilbert space $\Xi$.

\noindent So, at least formally, 
$$
(\call f)(x)=\frac{1}{2}(Tr BB^* \nabla^2 f)(x)+\<Ax,\nabla f(x)\>+\<B G(x),\nabla f(x)\>.
$$
The aim of this paper is to consider the case
where $\psi$ has superquadratic growth with respect to $\nabla v(t,x)B$, and both $\psi$
and $\phi$ are not bounded with respect to $x$. We consider the case of $\psi$ and $\phi$
differentiable, as well as the case of $\psi$ and $\phi$ only Lipschitz continuous: in this less regular case, in order to solve the
Kolmogorov equation \ref{Kolmo intro} we have to assume some regularizing property on the transition semigroup
$P_{t,\tau}$, namely
for every bounded and continuous real function $f$ on  $H$, for every $0<t<\tau\leq T$, the function $P_{t,\tau}\left[
f\right]  $ is differentiable
with respect to $x$ in directions in $\operatorname{Im}(B)$ (see also section \ref{sezionesde} for a detailed definition of the directional derivative $\nabla^B$),
and, $\forall x\in H,\,\xi\in \Xi$, 
\begin{equation}
\left|  \nabla P_{t,\tau}\left[  f\right]
\left(  x\right) B \xi\right|  \leq\frac{c}{\left(\tau-t\right)  ^{\alpha}%
}\left\|  f\right\|  _{\infty}\left|  \xi\right|. \label{ipotesi H intro}%
\end{equation}
We apply the results on equation (\ref{Kolmo intro}) to a stochastic optimal control problem.
Let us consider the controlled equation
\begin{equation}
\left\{
\begin{array}
[c]{l}%
dX^{u}_\tau  =\left[  AX^{u}_\tau +BG(X^u_\tau)+B R( u_\tau)  
 \right]  d\tau+BdW_\tau ,\text{ \ \ \ }\tau\in\left[  t,T\right] \\
X^{u}_t  =x.
\end{array}
\right.  \label{sdecontrol intro}%
\end{equation}
where $u$ is the control, taking values in a closed set $K$ of a normed space $U$.
Beside equation (\ref{sdecontrol intro}) we define the cost
\begin{equation*}
J\left(  t,x,u\right)  =\mathbb{E}\int_{t}^{T}[
\bar g\left(s,X^{u}_s\right)+g(u_s) ]ds+\mathbb{E}\phi\left(X^{u}_T\right),
\end{equation*}
for real functions $\mathbb{\phi}$, $\bar g$ and $g$. The control problem is to
minimize this functional $J$ over all admissible controls $u$.
We notice that we can treat a control problem with unbounded controls,
and we require weak coercivity on the cost $J$.
Indeed, we assume that, for $1<q\leq 2$, we have
\[
0\leq g(u)\leq c(1+\vert u \vert )^q
\qquad\text{and}\quad g(u)\geq C \vert u \vert ^q 
\qquad \text{for every }u \in K \textrm{ such that } \vert u\vert \geq R,
\]
so that the Hamiltonian function
\begin{equation*}
 \psi(t,x,z):=\bar g(t,x)+h(z):=\bar g(t,x)+\inf_{u\in K}\left\{  g\left(u\right)
+zR(u)\right\}
\end{equation*}
has quadratic or superquadratic growth, with respect to $z$, of order $p\geq 2$, the conjugate
exponent of $q$, if $q\leq 2$.

Second order differential equations on Hilbert spaces have been extensively studied (see e.g. the monograph \cite{DP3})
and one of the main motivations for this study
in the non linear case is the connection with control theory: in many cases the value
function of a finite horizon stochastic optimal control problem is solution to such a partial differential
equation. To study mild solutions of semilinear Kolmogorov equations
(\ref{Kolmo intro}) with $\psi$ Lipschitz continuous there are two main approaches in the literature: an
analytic approach and a purely probabilistic approach.
In the first direction we mention the paper \cite{Go1}, where the main
assumption is the strong Feller property for the transition semigroup $P_t $.

\noindent The purely probabilistic approach is based on backward stochastic
differential equations (BSDEs in the following). No regularizing assumption on the transition semigroup is imposed,
on the contrary $\psi$ and $\phi$ are assumed differentiable and $\psi$
is assumed to be Lipschitz continuous with respect to $y$ and $z$. In this direction
we refer to the paper \cite{fute}, which
is the infinite dimensional extension of results in \cite{PaPe}.

As far as we know, locally Lipschitz semilinear Kolmogorov equations with locally Lipschitz
Hamiltonian functions have been first treated in \cite{Go2}: $\psi$ is assumed to be locally Lipschitz continuous
with respect to $z$, and $\phi$ is taken Lipschitz continuous; both $\psi$ and $\phi$ are assumed to be bounded with respect to $x$.
The results in \cite{Go2} are achieved by means of a detailed study
on weakly continuous semigroups, and making the assumption
that the transition semigroup $P_t$ is strong Feller.

\noindent In \cite{BriFu} infinite dimensional
Hamilton Jacobi Bellman equations with Hamiltonian quadratic with respect to $z$
are solved in mild sense
by means of BSDEs: the generator $\call$ is related to
a more general Markov process $X$ than the one considered here in
(\ref{ornstein-pert intro}), and no regularizing assumptions on the coefficient are made,
but only the case of final condition $\phi$ G\^ateaux differentiable and bounded is treated.

\noindent In \cite{Mas3} infinite dimensional
Hamilton Jacobi Bellman equations with superquadratic Hamiltonian functions are considered, with bounded final condition,
also in the case of Lipschitz continuous coefficients, by requiring on the transition
semigroup the regularizing property we also mention in (\ref{ipotesi H intro}). Moreover in some special cases
the quadratic case is taken into account, with final condition only bounded and continuous.

\noindent In the present paper we improve the results both of \cite{Go2} and of \cite{BriFu}: we
are able to treat superquadratic Hamiltonian functions with an unbounded final condition, without requiring any regularizing properties
on the transition semigroup if the coefficients are G\^ateaux differentiable. We are also able
to take into account the case of Lipschitz continuous coefficients by requiring on $P_t$
the regularizing property already mentioned in (\ref{ipotesi H intro}), which is weaker than the strong Feller property
assumed in \cite{Go2}, and which is the same regularizing property considered in \cite{Mas3}.

Coming into the details of the techniques, in order to prove existence 
and uniqueness of a mild solution $v$ of equation (\ref{Kolmo intro}),
we use the fact that $v$ can be represented in terms of the solution of a suitable decoupled
forward-backward system (FBSDE in the following):
\begin{equation}\label{fbsde intro}
    \left\{\begin{array}{l}\dis
dX_\tau^{t,x}  =AX_\tau^{t,x} d\tau+B G(\tau,X_\tau^{t,x}) d\tau+BdW_\tau
,\text{ \ \ \ }\tau\in\left[  t,T\right], \\
X_{\tau}^{t,x} =x, \text{ \ \ \ }\tau\in\left[  0,t\right],\\ \dis
 dY_\tau^{t,x}=-\psi(\tau, X^{t,x}_\tau,Y_\tau^{t,x},Z_\tau^{t,x})\;d\tau+Z^{t,x}_\tau\;dW_\tau,
 \qquad \tau\in [0,T],
  \\\dis
  Y_T^{t,x}=\phi(X_T^{t,x}).
\end{array}\right.
\end{equation}
It is well known, see again \cite{PaPe} for the finite dimensional case and \cite{fute}
for the generalization to the infinite dimensional case, that $v(t,x)=Y_t^{t,x}$
when $\psi$ is Lipschitz continuous and all the data are differentiable. In \cite{BriFu} it is shown that
this identification holds true also when $\psi$ is quadratic and all the data differentiable, and in \cite{Mas3} it is further
extended, in the case of final datum bounded, to $\psi$ superquadratic and data not necessarily differentiable.
In this paper we go on extending this identification also in the case of final datum $\phi$
and Hamiltonian $\psi$ unbounded with polynomial growth with respect to $x$.

\noindent By the identification $v(t,x)=Y_t^{t,x}$, we achieve estimates on $v$ by studying the FBSDE
(\ref{fbsde intro}): we start from the results in \cite{Ri1}, and we extend them to the case
when the process $X^{t,x}$ solution of the forward equation in the FBSDE (\ref{fbsde intro})
takes values in an infinite dimensional Hilbert space $H$. The fundamental estimate we get is
 \begin{equation*}
\vert Z_\tau^{t,x}\vert\leq C(1+ \vert X_\tau^{t,x}\vert^r), \quad \forall \tau\in[0,T].
 \end{equation*}
where $r+1$ is the growth of $\psi $ with respect to $z$. From this estimate we deduce, in the case of differentiable
coefficients,
\begin{equation*}
\vert \nabla^Bv(t,x)\vert\leq C(1+ \vert x\vert^r), \quad \forall x\in H,\, t\in[0,T].
\end{equation*}
This is the fundamental tool to solve the HJB equation (\ref{Kolmo intro}) with differentiable coefficients.

\noindent To face the case of Lipschitz continuous coefficients,
and prove existence and uniqueness of a mild solution of equation (\ref{Kolmo intro}), assumption (\ref{ipotesi H intro})
on the transition semigroup $P_{t,\tau}$ is needed.
This condition is satisfied, among many other cases (see \cite{Mas}), by a stochastic wave equation on
the interval $[0,1]$
\begin{equation}\left\{
\begin{array}
[c]{l}%
\frac{\partial^{2}}{\partial\tau^{2}}y_\tau\left(\xi\right)  =\frac
{\partial^{2}}{\partial\xi^{2}}y_\tau\left(\xi\right)+f\left(\xi,y_\tau(\xi)\right)  
  +\dot{W}_\tau\left(\xi\right), \\
y_\tau\left(0\right)  =y_\tau\left(1\right)  =0,\\
y_t\left(\xi\right)  =x_{0}\left(  \xi\right)  ,\\
\frac{\partial y_\tau}{\partial\tau}\left( \xi\right)\mid_{\tau=t}  =x_{1}\left(  \xi\right).
\end{array}
\right.  \label{waveequation intro}
\end{equation}
Equation (\ref{waveequation intro}) can be reformulated in $H=L^{2}\left(  \left[  0,1\right]  \right)
\oplus H^{-1}\left(  \left[  0,1\right]  \right) $
as a perturbed Ornstein-Uhlenbeck process like (\ref{ornstein-pert intro}).

The paper is organized as follows:
in section \ref{sezionesde} we state notations and we recall some preliminary results on the
perturbed Ornstein-Uhlenbeck process, in section \ref{estimatesBSDE} we prove some fundamental estimates on the solution of the FBSDE,
in section \ref{sec-diffle} we study differentiability of the Markovian BSDE (\ref{fbsde intro})
when all the coefficients are differentiable. Thanks to these results we are able to solve Kolmogorov equation
(\ref{Kolmo intro}) with differentiable coefficients, see section \ref{sezionePDEdiffle}.
In section \ref{sezionePDElip} we turn to only locally Lipschitz continuous coefficients. Finally in section \ref{applic contr}
we apply the results to a finite horizon optimal control problem, and in \ref{sez-contr-wave} we present
the special case of a controlled wave equation.

\section{Notations and preliminary results}
\label{sezionesde}

\subsection{Notations}

In this paper we denote by $H$ and $\Xi$ some real and separable Hilbert spaces, and by $(W_t)_{t\geq 0}$
a cylindrical Wiener process in $\Xi$, defined on a complete probability
space $(\Omega, \calf, \P)$. For $t\geq 0$, let $\calf_t$ denote the $\sigma$-algebra
generated by $(W_s,\, s\leq t)$ and augmented with the $\P$-null sets of $\calf$. The notation
$\E_t$ stands for the conditional expectation given $\calf_t$.

\noindent Given a real and separable Hilbert space $K$, (eventually $K=\R^m$), we denote further
\begin{itemize}
 \item $\cals^p(K)$, or $\cals^p$ where no confusion is possible, the space of all adapted and c\`adl\`ag processes
$(Y_t)_{t\in[0,T]}$ with values in $K$, normed by 
$\Vert Y\Vert _{\cals^p}=\E[\sup_{t\in[0,T]}\vert Y_t\vert^p]^{1/p}$; $\cals^\infty(K)$,
or $\cals^\infty$ where no confusion is possible, the space of all bounded predictable processes.
\item $\calm^p(K)$, or $\calm^p$ where no confusion is possible, the space of all predictable processes
$(Z_t)_{t\in[0,T]}$ with values in $K$, normed by 
$\Vert Z\Vert _{\calm^p}=\E[(\int_0^T\vert Z_t\vert^2dt)^{p/2}]^{1/p}$.
\end{itemize}

We recall that a function $f : X \rightarrow V$ where $X$ and $V$ are two Banach spaces, has a directional derivative at point $x \in X$ in the direction $h \in X$ when
$$\nabla f(x;h) =\lim_{s \rightarrow 0} \frac{f(x+sh)-f(x)}{s}$$
exists. $f$ is said to be G\^ateaux differentiable at point $x$ if $\nabla f(x;h)$ exists for every $h$ and there exists an element of $L(X,V)$, denoted as $\nabla f(x)$ and called the G\^ateaux derivative, such that $\nabla f(x;h)=\nabla f(x) h$ for every $h \in X$. Let us introduce some notations.
\begin{itemize}
 \item $f: X \rightarrow V$ belongs to the class $\mathcal{G}^1(X;V)$ if it is continuous, G\^ateaux differentiable on $X$, and $\nabla f : X \rightarrow L(X,V)$ is strongly continuous.
 \item $f: X \times Y \rightarrow V$ belongs to the class $\mathcal{G}^{1,0}(X \times Y;V)$ if it is continuous, G\^ateaux differentiable with respect to its first variable $x \in X$, and $\nabla_x f : X\times Y \rightarrow L(X,V)$ is strongly continuous.
\end{itemize}
When $f$ depends on additional arguments, the previous definitions have obvious generalizations.

We briefly introduce the notion of $B$-differentiability,
for further details see e.g. \cite{Mas}. We recall that for a continuous function
$f:H\rightarrow\mathbb{R}$ the $B$-directional derivative $\nabla^B$ at a
point $x\in H$ in direction$\ \xi\in H$ is defined as follows:%
\[
\nabla^{B}f\left(  x;\xi\right)  =\lim_{s\rightarrow0}\frac{f\left(
x+sB\xi\right)  -f\left(  x\right)  }{s},\text{ }s\in\mathbb{R}\text{.}%
\]
A continuous function $f$ is $B$-G\^ateaux differentiable at a point $x\in H$ if
$f$ admits the $B$-directional derivative $\nabla^Bf\left(  x;\xi\right)  $
in every directions $\xi\in \Xi$ and there exists a functional, the
$B-$gradient $\nabla^Bf\left(  x\right)  \in\Xi^{\ast}$ such that
$\nabla^Bf\left(  x;\xi\right)  =\nabla^Bf\left(  x\right)  \xi$.

Finally, $C$ will denote, as usual, a constant that may change its value from line to line.

\subsection{The forward equation}

We consider a perturbed Ornstein-Uhlenbeck process in 
$H$, that is a Markov process $X$ (also denoted $X^{t,x}$ to stress the dependence
on the initial conditions) solution to equation%
\begin{equation}
\left\{
\begin{array}
[c]{l}%
dX_\tau  =AX_\tau d\tau+F(\tau,X_\tau) d\tau+BdW_\tau
,\text{ \ \ \ }\tau\in\left[  t,T\right], \\
X_{\tau} =x,\text{ \ \ \ }\tau\in\left[  0,t\right],
\end{array}
\right.  \label{ornstein-pert}%
\end{equation}
where $A$ is the generator of a strongly continuous semigroup in $H$,
$B$ is a linear bounded operator from $\Xi$ to $H$ and 
$F$ is a map from $[0,T]\times H$ with values in $H$.
\noindent We define the positive and symmetric operator%
\[
Q_{\sigma}=\int_{0}^{\sigma}e^{sA}BB^{\ast}e^{sA^{\ast}}ds.
\]
Throughout the paper we assume the following.

\begin{hypothesis}
\label{ip su AB}

\begin{enumerate}
\item  The linear operator $A$ is the generator of a strongly continuous
semigroup $\left(  e^{t A},t\geq0\right)  $ in the Hilbert space $H.$
It is well known that there exist $N>0$ and $\omega\in\mathbb{R}$ such that
$\left\Vert e^{tA}\right\Vert _{L\left(  H,H\right)  }\leq Ne^{\omega t}$, for
all $t\geq0$. In the following, we always consider $N\geq1$ and $\omega\geq 0$.
\item $F:[0,T]\times H\longrightarrow H$ is continuous and $\forall\, \tau \in [0,T]$,
$F(\tau,.)$ is Lipschitz continuous and belongs to $\mathcal{G}^1(H,H)$: $\forall\, \tau \in [0,T]$
and $\forall x,x' \in H$
\[
\vert F(\tau,0)\vert\leq C; \qquad \vert F(\tau,x)-F(\tau,x')\vert \leq K_F\vert x-x'\vert.
\]
As a consequence, $\vert\nabla_x F(\tau,x)\vert\leq K_F$, $\forall \tau\in[0,T],\,
\forall x\in H$.
\item $B$ is a bounded linear operator from $\Xi$ to $H$ and $Q_{\sigma}$ is
of trace class for every $\sigma\geq0$.
\end{enumerate}
\end{hypothesis}
We notice that the differentiability assumption on $F$ will be used to prove differentiability of the process $X$ with respect to the initial datum $x$,
and it is not necessary to prove existence of a solution to equation \ref{ornstein-pert},
which is a standard result collected in the following proposition.
\begin{proposition}
\label{esistenza sde}Under Hypothesis \ref{ip su AB}, the forward equation in (\ref{ornstein-pert})
admits a unique continuous mild solution. Moreover 
\[\mathbb{E}\left[\sup_{\tau
\in\left[  0,T\right]  }\left|  X_\tau^{t,x}  \right|  ^{p}\right]%
<C_{p}\left(  1+\left|  x\right|  \right)  ^{p},\] for every $p\in\left(
0,\infty\right)  ,$ and some constant $C_{p}>0$.
\end{proposition}

\noindent If $F=0$, $X^{t,x}$ is an Ornstein-Uhlenbeck process and it
is clearly time-homogeneous, and for
$0\leq t\leq \tau \leq T$ we denote by $P_{\tau-t}=P_{t,\tau}$
its transition semigroup, where for every bounded and continuous function $\phi:H\rightarrow\R$
\[
 P_{t,\tau}[\phi](x)=\E\phi(X_\tau^{t,x}).
\]
It is well known that the
Ornstein-Uhlenbeck semigroup can be represented as
\[
P_{\sigma}\left[  \phi\right]  \left(  x\right)  :=\int_{H}\phi\left(  y\right)
\mathcal{N}\left(  e^{\sigma A}x,Q_{\sigma}\right)  \left(  dy\right), \quad \sigma >0  ,
\]
where $\mathcal{N}\left(  e^{\sigma A}x,Q_{\sigma}\right)  \left(  dy\right)  $ denotes a
Gaussian measure with mean $e^{\sigma A}x,$ and covariance operator $Q_{\sigma}$.

\section{Some estimates on (super)quadratic BSDEs in infinite dimensional Markovian framework}
\label{estimatesBSDE}
In this section we consider the following BSDE
\begin{equation}\label{bsde}
    \left\{\begin{array}{l}\dis
 dY_\tau^{t,x}=-\psi(\tau, X^{t,x}_\tau,Y_\tau^{t,x},Z_\tau^{t,x})\;d\tau+Z^{t,x}_\tau\;dW_\tau,
 \qquad \tau\in [0,T],
  \\\dis
  Y_T^{t,x}=\phi(X_T^{t,x}),
\end{array}\right.
\end{equation}
where $X^{t,x}$ is a perturbed Ornstein-Uhlenbeck process solution of equation (\ref{ornstein-pert}).
We call it also BSDE in Markovian framework and we note that in this paper
$X$ is an infinite dimensional Markov process.
Under suitable assumptions on the coefficients
 $\psi:[0,T]\times H \times \mathbb{R} \times \Xi^{\ast}
\rightarrow\mathbb{R}$ 
and $\mathbb{\phi}:H\rightarrow\mathbb{R}$
we will look for a solution consisting of a pair of predictable processes,
taking values in $\mathbb{R}\times \Xi$, such that $Y$ has
continuous paths and
\[
\|\left( Y,Z\right)\|_{\cals^2\times\calm^2}<\infty.
\]

We make the following assumptions on the generator $\psi$ and on the final
datum $\phi$ in the backward equation (\ref{bsde}).
\begin{hypothesis}
\label{ip-psi-phi}The maps $\phi:H
\rightarrow\mathbb{R}$, $\psi:[0,T]\times H\times \R\times \Xi^{\ast}\rightarrow\mathbb{R}$ are continuous
and there exist constants $l\geq1,\;0\leq r<\frac{1}{l},\;\alpha\geq 0,\;\beta\geq 0\; ,\gamma\geq0$ and
$K_{\psi_y}>0$ such that
\begin{enumerate}
\item for all $(t,x,y,y',z)\in [0,T]\times H\times \R\times \R\times \Xi^{\ast}$,
\[
 \vert \psi(t,x,y,z)-\psi(t,x,y',z)\vert \leq K_{\psi_y}\vert y -y'\vert;
\]
\item for all $(t,x,y,z,z')\in [0,T]\times H\times \R\times \Xi^{\ast}\times \Xi^{\ast}$
\[
 \vert \psi(t,x,y,z)-\psi(t,x,y,z')\vert \leq 
\left(C+\frac{\gamma}{2}\vert z\vert^l+\frac{\gamma}{2}\vert z'\vert^l\right)\vert z -z'\vert;
\]
\item for all $(t,x,x',y,z)\in [0,T]\times H\times H\times \R\times \Xi^{\ast}$
\[
 \vert \psi(t,x,y,z)-\psi(t,x',y,z)\vert \leq 
\left(C+\frac{\beta}{2}\vert x\vert^r+\frac{\beta}{2}\vert x'\vert^r\right)\vert x -x'\vert;
\]
\[
 \vert \phi(x)-\phi(x')\vert \leq 
\left(C+\frac{\alpha}{2}\vert x\vert^r+\frac{\alpha}{2}\vert x'\vert^r\right)\vert x -x'\vert.
\]
\end{enumerate}
\end{hypothesis}
Notice that in previous assumptions the quadratic case corresponds to $l=1$
 and the superquadratic case to $l>1$.
Before proving an existence and uniqueness result for the BSDE (\ref{bsde}), we prove the following lemma.
\begin{lemma}\label{lemma-exist-bsde}
 Assume that Hypothesis \ref{ip su AB} holds true.
Moreover, we assume on the final datum and on the generator of the BSDE (\ref{bsde}) that:
\begin{itemize}
 \item $\phi: H \rightarrow H$ is a Lipschitz continuous function with Lipschitz constant given by $K_\phi$;
\item $\psi : [0,T]\times H \times \mathbb{R} \times \Xi^{\ast}$ is a continuous function and there exist constants
$K_{\psi_x},\,K_{\psi_y},\,K_{\psi_z}$ such that
$\forall \tau \in [0,T],\,x,x'\in H,y,y'\in\R,z,z'\in\Xi^{\ast}$
\begin{align*}
 &\vert \psi(\tau,x,y,z)-\psi(\tau,x',y,z)\vert \leq K_{\psi_x}\vert x -x'\vert;\\ \nonumber
 &\vert \psi(\tau,x,y,z)-\psi(\tau,x,y',z)\vert \leq K_{\psi_y}\vert y -y'\vert;\\ \nonumber
 &\vert \psi(\tau,x,y,z)-\psi(\tau,x,y,z')\vert \leq K_{\psi_z}(1+\varphi(\vert z\vert)+
\varphi(\vert z'\vert))\vert z -z'\vert;\\ \nonumber
\end{align*}
\end{itemize}
where $\varphi:\R^+\rightarrow\R^+$ is a non decreasing function.
Then the BSDE (\ref{bsde}) admits a unique solution $(Y,Z)\in\cals^2\times\calm^2$
such that
\begin{equation*}
 \vert Z_\tau^{t,x}\vert \leq C,
\end{equation*}
where $C$ is a constant depending on $A,\,F,\,B,\,K_\phi,\,K_{\psi_x},\,K_{\psi_y}$ and $T$.
\end{lemma}

\noindent {\bf Proof.}
We use a classical truncation argument:
we set $\psi_M=\psi (\cdot,\cdot,\cdot,\rho_M(\cdot))$,
where $\rho_M$ is a smooth modification of the projection on the centered ball of radius $M$
such that $\vert\rho_M\vert\leq M$, $\vert\nabla\rho_M\vert\leq 1$ and $\rho_M(x)=x$ when
$\vert x\vert\leq M-1$. In particular $\psi_M$ is also Lipschitz continuous with respect to $z$.
Now assuming first that $\phi$ and $\psi$ are differentiable with respect to
$x$, $y$ and $z$, it turns out that $\psi_M$ is also differentiable
with respect to $x$, $y$ and $z$. So we can differentiate
the BSDE
\begin{equation*}
    \left\{\begin{array}{l}\dis
 dY_\tau^{M,t,x}=-\psi_M(\tau, X^{t,x}_\tau,Y_\tau^{M,t,x},Z_\tau^{M,t,x})\;d\tau+Z^{M,t,x}_\tau\;dW_\tau,
 \qquad \tau\in [0,T],
  \\\dis
  Y_T^{M,t,x}=\phi(X_T^{t,x}),
\end{array}\right.
\end{equation*}
with respect to the initial condition $x$ in the forward equation (\ref{ornstein-pert}).
We get
\begin{equation*}
    \left\{\begin{array}{l}\dis
 d \nabla_xY_\tau^{M,t,x}=-\nabla_x\psi_M(\tau, X^{t,x}_\tau,Y_\tau^{M,t,x},Z_\tau^{M,t,x})\nabla_xX^{t,x}_\tau\;d\tau
\\ \dis
\qquad \qquad\qquad-\nabla_y\psi_M(\tau, X^{t,x}_\tau,Y_\tau^{M,t,x},Z_\tau^{M,t,x})\nabla_xY^{M,t,x}_\tau\;d\tau
\\ \dis
\qquad \qquad\qquad-\nabla_z\psi_M(\tau, X^{t,x}_\tau,Y_\tau^{M,t,x},Z_\tau^{M,t,x})\nabla_xZ^{M,t,x}_\tau\;d\tau
+\nabla_xZ^{M,t,x}_\tau\;dW_\tau,
 \qquad \tau\in [0,T],
  \\\dis
  \nabla_x Y_T^{M}=\nabla\phi(X_T^{t,x})\nabla X^{t,x}_T.
\end{array}\right.
\end{equation*}
Since $\nabla_z\psi_M(\tau, X^{t,x}_\tau,Y_\tau^{M,t,x},Z_\tau^{M,t,x})$ is bounded, 
we can apply Girsanov's theorem: there exists
a probability measure $\Q^M$,
equivalent to the original one $\P$, such that
$$\tilde{W}_t:=W_t-\int_0^t \nabla_z\psi_M(s,X_s,Y^{M,t,x}_s,Z^{M,t,x}_s)ds$$ is a Wiener process under $\Q^M$.
We obtain
\begin{align*}
 \nabla Y^{M,t,x}_\tau =&\E_\tau^{\Q^M}\bigg[ e^{\int_\tau^T
\nabla_y\psi_M(u,X_u^{t,x},Y^{M,t,x}_u,Z^{M,t,x}_u)du}\nabla\phi(X_T^{t,x})\nabla X_T^{t,x}\\ \nonumber
&+\int_\tau^Te^{\int_\tau^s
\nabla_y\psi_M(u,X_u^{t,x},Y^{M,t,x}_u,Z^{M,t,x}_u)du}\
\nabla_x\psi_M(s,X_s^{t,x},Y^{M,t,x}_s,Z^{M,t,x}_s)\nabla X_s^{t,x}ds\bigg],
\end{align*}
from which we deduce a bound, uniform with respect to $x$, $t$ and $\tau$, for $ \nabla Y^{M,t,x}_\tau$,
and consequently for $ \nabla Y^{M,t,x}_\tau \, B$:
\begin{equation}\label{stimanablaYM}
\vert  \nabla Y^{M,t,x}_\tau  B\vert\leq C,
\end{equation}
with $C$ a constant which does not depend on $x$, $t$ and $\tau$.
By the Markov property (see e.g. part 5 in \cite{fute}), we have
$$Z_{\tau}^{M,t,x} = Z_{\tau}^{M,\tau,X_{\tau}^{t,x}}=Z_{\tau}^{M,\tau,y}|_{y=X_{\tau}^{t,x}}.$$
In \cite{fute}, a standard result on BSDEs with Lipschitz generator in infinite dimensional framework gives us also that $$\nabla Y_{\tau}^{M,\tau,y}B=Z_{\tau}^{M,\tau,y}.$$
Finally, by using estimate (\ref{stimanablaYM}), we obtain
$$\abs{Z_{\tau}^{M,t,x}}=\abs{Z_{\tau}^{M,\tau,y}}_{y=X_{\tau}^{t,x}}=\abs{\nabla Y_{\tau}^{M,\tau,y} B}_{y=X_{\tau}^{t,x}} \leqslant C,$$
with $C$ a constant that does not depend on $M$. So, for $M$ large enough we have $\rho_M(Z^{M,t,x})=Z^{M,t,x}$ and $(Y^{M,t,x},Z^{M,t,x})$ becomes a solution of the initial BSDE (\ref{bsde}). The uniqueness comes from the classical uniqueness result for Lipschitz BSDEs.

Notice that, unlike in finite dimensions, we cannot consider, for any
$s\in[t,T]$ $\left(\nabla X_s^{t,x}\right)^{-1} $, unless $A$
is the generator of a group, while in the present paper we consider
with more generality that $A$ is the generator of a semigroup.

\qed

\noindent Now we are ready to prove an existence and uniqueness result for the BSDE (\ref{bsde}),
together with an estimate on $Z$,
when $X^{t,x}$ is an Ornstein-Uhlenbeck process, that is to say $F=0$ in equation (\ref{ornstein-pert}).
We essentially follow the proof of Proposition 2.2 in \cite{Ri1}, with suitable differences
due to the infinite dimensional setting.
\begin{proposition}\label{prop-exist-bsde}
Assume that Hypotheses \ref{ip su AB}, with $F=0$, and \ref{ip-psi-phi} hold true.
Then there exists  a solution $(Y^{t,x},Z^{t,x})$ of the Markovian BSDE (\ref{bsde}) such that
$(Y^{t,x},Z^{t,x})\in \cals^2\times\calm^2$ and
 \begin{equation}\label{estimate-Z-r}
\vert Z^{t,x}_\tau\vert\leq C(1+ \vert X^{t,x}_\tau\vert^r), \quad \forall \tau\in[0,T].
 \end{equation}
Moreover this solution is unique amongst solutions such that
\begin{itemize}
 \item $Y^{t,x}\in\cals^2$;
\item there exists $\eta>0$ such that 
\[
 \E\left[ e^{(\frac{1}{2}+\eta)\frac{\gamma^2}{4}\int_0^T\vert Z^{t,x}_s\vert^{2l}ds}\right]<+\infty.
\]
\end{itemize}
\end{proposition}
\noindent {\bf Proof.}  We remark that if there exists a solution $(Y^{t,x},Z^{t,x})$ such that
 \[
  \vert Z^{t,x}_\tau\vert\leq C\left(1+\vert X^{t,x}_\tau\vert^r\right), \quad \forall \tau\in[0,T],
 \]
then, $\forall c>0$,
\begin{align*}
 \E&\left[ e^{c\int_0^T\vert Z^{t,x}_s\vert^{2l}ds}\right]
\leq C\E\left[ e^{C\int_t^T\vert X_s^{t,x}\vert^{2lr}\,ds}\right]=C\E\left[ e^{(T-t)\times\frac{1}{T-t}C\int_t^T\vert X_s^{t,x} \vert^{2lr}\,ds}\right]\\ \nonumber
&\leq C\E\left[ \frac{1}{T-t}\int_t^Te^{C(T-t)\vert X_s^{t,x} \vert^{2lr}}ds\right]\leq C \frac{1}{T-t}\int_t^T\E\left[e^{CT\vert X_s^{t,x} \vert^{2lr}}\right]ds<+\infty,\nonumber
\end{align*}
where we have used Jensen inequality. The last bound follows from inequality $2lr<2$
(see assumptions in Hypothesis \ref{ip-psi-phi})
and the fact that $X$ is an Ornstein-Uhlenbeck process, and so in particular a Gaussian random variable:
$X_s^{t,x}\sim\caln(e^{(s-t)A}x,Q_{s-t})$.

\noindent Now uniqueness follows as in the proof of proposition 2.2 in \cite{Ri1}.

\noindent For what concerns existence, following again \cite{Ri1},
we approximate the Markovian BSDE (\ref{bsde}) by a truncation argument, namely we consider
$(Y^{M,t,x},Z^{M,t,x})$ solution of the following BSDE
\begin{equation}\label{bsde-M}
    \left\{\begin{array}{l}\dis
 dY_\tau^{M,t,x}=-\psi_M(\tau, X^{t,x}_\tau,Y_\tau^{M,t,x},Z_\tau^{M,t,x})\;d\tau+Z^{M,t,x}_\tau\;dW_\tau,
 \qquad \tau\in [0,T],
  \\\dis
  Y_T^{M,t,x}=\phi_M(X_T^{t,x}),
\end{array}\right.
\end{equation}
where $\phi_M=\phi\circ\rho_M$, $\psi_M=\psi (\cdot,\rho_M(\cdot),\cdot,\cdot)$,
and $\rho_M$ is a smooth modification of the projection on the centered ball of radius $M$
such that $\vert\rho_M\vert\leq M$, $\vert\nabla\rho_M\vert\leq 1$ and $\rho_M(x)=x$ when
$\vert x\vert\leq M-1$. So $g_M$ and $\psi_M$ are Lipschitz and bounded functions
with respect to $x$.
By Lemma \ref{lemma-exist-bsde}, we get that there exists a unique solution
$(Y^{M,t,x},Z^{M,t,x})$ to the BSDE (\ref{bsde-M}) such that $\abs{Z^{M,t,x}}\leq A_0$ with $A_0$ a constant that depends on $M$. As a consequence, $\psi_M$ is a Lipschitz function
with respect to $z$ and so classical results on BSDEs apply.
Next assume for a moment the following lemma, whose proof is
similar to the proof of lemma 2.4 in \cite{Ri1}, 
\begin{lemma}\label{lemma-iterate}
 Under assumptions of Proposition \ref{prop-exist-bsde}, we have
\begin{equation*}
\vert Z^{M,t,x}\vert\leq A_n+B_n\vert X^{t,x}\vert^r, 
\end{equation*}
with $(A_n,B_n)_{n\in\N}$ defined by recursion: $B_0=0$, $A_0$ defined before,
\begin{align}\label{relA_nB_n}
&A_{n+1}= C(1+A_n^{lr}),\\ \nonumber
&B_{n+1}=C,
\end{align}
where $C$ is a constant that does not depend on $M$.
\end{lemma}
Notice that relation (\ref{relA_nB_n}) for $A_n$ is a contraction, so
its limit exists, we denote it by $A_\infty$, and it does not depend on $M$,
so
\[
\vert Z^{M,t,x}\vert\leq A_\infty+C\vert X^{t,x}\vert^r.
\]
Now it remains to show that $(Y^{M,t,x},Z^{M,t,x})_{M \in \mathbb{N}}$ is a Cauchy sequence that tends to a limit $(Y^{t,x},Z^{t,x})$ solution of the BSDE (\ref{bsde}).
This part of the proof goes on like in \cite{Ri1}, Proposition 2.2.
\qed

\noindent {\bf Proof of Lemma \ref{lemma-iterate}.} The proof is similar to the
proof of Lemma 2.4 in \cite{Ri1}, and we give it for the reader convenience and to give references for the infinite dimensional setting.

\noindent We start by considering $\phi$ G\^ateaux differentiable and
$\psi$ G\^ateaux differentiable with respect to $x,\,y$ and $z$.
As in \cite{Ri1}, the proof is given by recursion: for $n=0$, by
lemma \ref{lemma-exist-bsde}, the result is true,
let us suppose that it is true for some $n\in\N$ and let us show that it is still true for $n+1$.
We get that $X^{t,x}$, $(Y^{M,t,x},Z^{M,t,x} )$
are differentiable. Arguing as in \cite{Ri1}, since
\[
 \vert\nabla_z\psi_M(s,X_s^{t,x},Y^{M,t,x}_s,Z^{M,t,x}_s)\vert\leq C(1+\vert Z^{M,t,x}_s\vert^l)\leq C_M,
\]
by the Girsanov theorem there exists a probability measure $\Q^M$,
equivalent to the original one $\P$, such that 
$\tilde{W}_{\tau}:=W_{\tau}-\int_0^{\tau} \nabla_z\psi_M(s,X_s^{t,x},Y^{M,t,x}_s,Z^{M,t,x}_s)ds$ is a Wiener process under $\Q^M$.
We obtain
\begin{align*}
 \nabla Y^{M,t,x}_\tau&=\E_\tau^{\Q^M}\left[ e^{\int_\tau^T
\nabla_y\psi_M(u,X_u^{t,x},Y^{M,t,x}_u,Z^{M,t,x}_u)du}\nabla\phi_M(X_T^{t,x})\nabla X_T^{t,x}\right.\\ \nonumber
&\left.+\int_\tau^Te^{\int_\tau^T
\nabla_y\psi_M(u,X_u^{t,x},Y^{M,t,x}_u,Z^{M,t,x}_u)du}\
\nabla_x\psi_M(s,X_s^{t,x},Y^{M,t,x}_s,Z^{M,t,x}_s)\nabla X_s^{t,x}ds\right],
\end{align*}
and, by using assumptions 3.1 and the fact that $\nabla X^{t,x}$ is bounded,
$$\abs{\nabla Y^{M,t,x}_\tau B} \leq C+C \E_\tau^{\Q^M}
\left[ \vert X_T^{t,x}\vert^r+ \int_{\tau}^T\vert X_s^{t,x}\vert^r  ds \right]. $$
Once again, the Markov property and standard results on BSDEs with Lipschitz generator in infinite dimension (see e.g. \cite{fute}) give us
\begin{eqnarray} 
\nonumber
 \abs{Z_{\tau}^{M,t,x}}&=&\abs{Z_{\tau}^{M,\tau,x'}}_{x'=X_{\tau}^{t,x}}=\abs{\nabla Y_{\tau}^{M,\tau,x'} B}_{x'=X_{\tau}^{t,x}}\\
\nonumber &\leq& \left.\left(  C+C \E_\tau^{\Q^M}
\left[ \vert X_T^{t,x'}\vert^r+ \int_t^T\vert X_s^{t,x'}\vert^r  ds \right]\right)\right|_{x'=X_{\tau}^{t,x}}\\
 &\leq&C+C \E_\tau^{\Q^M}
\left[ \vert X_T^{t,X_{\tau}^{t,x}}\vert^r+ \int_t^T\vert X_s^{t,X_{\tau}^{t,x}}\vert^r  ds \right]. \label{stima Z M}
\end{eqnarray}

Now we need to estimate $\E_\tau^{\Q^M}\vert X_s^{t,x}\vert^r$, for $s\in[\tau,T]$.
In $(\Omega,\calf,\Q^M)$, for $s\in[\tau,T]$, $X_s^{t,x}$ solves the following equation in mild form
\begin{equation}\label{ornstein-QM}
 X_s^{t,x}=e^{(s-\tau)A}X_\tau^{t,x}+\int_\tau^s e^{(s-r)A}B\,d\tilde{W_r}
+\int_\tau^s e^{(s-r)A}B \nabla_z\psi_M(r,X_r^{t,x},Y^{M,t,x}_r,Z^{M,t,x}_r)\,dr.
\end{equation}
Notice that $\E_\tau^{\Q^M}\abs{\int_\tau^s e^{(s-r)A}B\,d\tilde{W}_r}=
\E^{\Q^M}\abs{\int_\tau^s e^{(s-r)A}B\,d\tilde{W}_r}$,
and by Corollary 2.17 in \cite{DP1}, 
\begin{align*}
 \E^{\Q^M}\abs{\int_\tau^s e^{(s-r)A}B\,d\tilde{W_r}}\leq \left(\E^{\Q^M}\abs{\int_\tau^s e^{(s-r)A}B\,d\tilde{W_r}}^2\right)^{1/2}
\leq C \left(\int_{\tau}^s e^{(s-r)A}BB^*e^{(s-r)A^*}\,ds\right)^{1/2}<\infty.
\end{align*}
So, we have
\[
 \E^{\Q^M}\abs{\int_\tau^s e^{(s-r)A}B\,d\tilde{W_r}}\leq C
\]
where $C$ is a constant that depends on $A,\,B$.

Coming back to the estimate of $\E_\tau^{\Q^M}\vert X_s^{t,x}\vert$, we get, using the last inequality and (\ref{ornstein-QM}),
\begin{align*}
\E_\tau^{\Q^M}\vert X_s^{t,x}\vert&\leq
Ne^{(s-\tau)\omega}\vert X_\tau^{t,x}\vert +C
+N\E_\tau^{\Q^M}\int_\tau^s e^{(s-\sigma)\omega}\Vert B\Vert_{L(\Xi,H)}\gamma\vert Z^{M,t,x}_\sigma\vert \,d\sigma \\ \nonumber
&\leq
Ne^{(s-\tau)\omega}\vert X_\tau^{t,x}\vert +C
+N\E_\tau^{\Q^M}\int_\tau^s e^{(s-\sigma)\omega}\Vert B\Vert_{L(\Xi,H)}
\gamma(A_n+B_n\vert X^{t,x}_\sigma\vert^{r}\vert)^l \,d\sigma \\ \nonumber
&\leq C\vert X_\tau^{t,x}\vert +C+CA_n^l+C\E_\tau^{\Q^M}\int_\tau^s B_n^l\vert X^{t,x}_\sigma\vert^{rl} \,d\sigma.
\end{align*}
By applying Young inequality, with $1/p+1/q=1$ and $rlq=1$, we get
\[
C B_n^l\vert X^{t,x}_\sigma\vert^{rl}\leq
\frac{C^p B_n^{lp}}{p}+\frac{\vert X^{t,x}_\sigma\vert^{rlq}}{q}.
\]
Thus, we have
\[
\E_\tau^{\Q^M}\vert X_s^{t,x}\vert\leq C\vert X_\tau^{t,x}\vert +C+CA_n^l+C B_n^{lp}+C\E_\tau^{\Q^M}\int_\tau^s \vert X^{t,x}_\sigma\vert \,d\sigma.
\]
Finally, by Gr\"onwall lemma we get
\[
\E_\tau^{\Q^M}\vert X_s^{t,x}\vert\leq e^{C(T-\tau)}\left( C\vert X_\tau^{t,x}\vert +C+CA_n^l+C B_n^{lp}\right),
\]
and, since $r<1$,
\[
\E_\tau^{\Q^M}\vert X_s^{t,x}\vert^r\leq C\left(1+A_n^{lr}+B_n^{lrp}+\vert X_\tau^{t,x}\vert^r \right).
\]
Thus, by (\ref{stima Z M}) and the recursion assumption on $B_n$, we get
\begin{equation*}
\vert Z^{M,\tau,y}\vert \leq C\left(1+A_n^{lr}+B_n^{lrp}+\vert X_\tau^{t,x}\vert^r \right) \leq C\left(1+A_n^{lr}+\vert X_\tau^{t,x}\vert^r \right),
\end{equation*}
so we can take
\[
B_{n+1}=C
\]
and
$$ A_{n+1}=C\left( 1+A_n^{rl} \right),$$
so that $(B_n)_{n\in\N^*}$
is a constant sequence and $(A_n)_{n\in\N}$ satisfies the required recursion relation.
When $\phi$ and $\psi$ are not G\^ateaux differentiable, we can approximate
them by their inf-sup convolutions,
noting that since $\phi_M$ and $\psi_M$ are Lipschitz continuous,
also their inf-sup convolutions are.
\qed

Now we prove an analogous of Proposition \ref{prop-exist-bsde} when $X$
is a perturbed Ornstein-Uhlenbeck process.

\begin{proposition}\label{prop-exist-bsde-pert}
Assume that Hypotheses \ref{ip su AB} and \ref{ip-psi-phi} hold true.
Then there exists a solution $(Y^{t,x},Z^{t,x})$ of the Markovian BSDE (\ref{bsde}) such that
$(Y^{t,x},Z^{t,x})\in \cals^2\times\calm^2$ and
\begin{equation}
 \label{estimate-Z-r-Fneqzero}
 \vert Z^{t,x}_\tau\vert\leq C\left(1+ \vert X^{t,x}_\tau\vert^r \right), \quad \forall \tau\in[0,T].
\end{equation}
Moreover this solution is unique amongst solutions such that
\begin{itemize}
 \item $Y\in\cals^2$;
\item there exists $\eta>0$ such that 
\[
 \E\left[ e^{(\frac{1}{2}+\eta)\frac{\gamma^2}{4}\int_0^T\vert Z^{t,x}_s\vert^{2l}ds}\right]<+\infty.
\]
\end{itemize}
\end{proposition}
\noindent {\bf Proof.}  We only give the proof of the points where some
differences with the case of an Ornstein-Uhlenbeck process treated in Proposition \ref{prop-exist-bsde} arise.

\noindent As a first point, let us prove that, for all $p>0$ there exists a constant $C$ that does not depend on $(t,x)$ such that
\begin{equation}
\label{exponential moment sup X2}
\E\left[ e^{p\int_0^T\vert X_s^{t,x}\vert^{2lr}ds}\right]\leq Ce^{C\abs{x}^{2rl}}<+\infty
\end{equation}
where this time
$X^{t,x}$ satisfies, in mild form,
\[
 X_\tau^{t,x}  =e^{(\tau-t)A}x+\int_t^\tau e^{(\tau-s)A}F(s,X_s^{t,x}) ds+
\int_t^\tau e^{(\tau-s)A}BdW_s.
\]
To prove (\ref{exponential moment sup X2}) we denote the stochastic convolution by $W_A(\tau):=\dis\int_t^\tau e^{(\tau-s)A}BdW_s$ and we set
\[
 \Gamma_\tau^{t,x}:=X_\tau^{t,x}-W_A(\tau).
\]
The process $\Gamma$ satisfies the integral equation
\[
 \Gamma_\tau^{t,x}  =e^{(\tau-t)A}x+\int_t^\tau e^{(\tau-s)A}F(s,\Gamma_s^{t,x}+W_A(s)) ds.
\]
So
\begin{eqnarray*}
\vert \Gamma_\tau^{t,x}\vert &\leq& Ne^{\omega T}\vert x\vert+N K_F\int_t^\tau e^{\omega(\tau-s)}
\left(1+\vert \Gamma_s^{t,x}\vert+\vert W_A(s)\vert\right) ds\\
&\leq & C\left(1+\vert x\vert+\int_t^\tau \vert W_A(s)\vert ds\right)+C\int_t^\tau\vert \Gamma_s^{t,x}\vert ds.
\end{eqnarray*}
By a generalization of the Gronwall lemma in integral form we get
\[
\vert \Gamma_\tau^{t,x}\vert \leq C\left(1+\abs{x}+\int_t^\tau \vert W_A(s)\vert ds\right)+C\int_t^{\tau}\int_t^u\abs{W_A(s)} ds du,
\]
so,
\[
\vert \Gamma_s^{t,x}\vert^{2lr}
\leq C\left(1+\abs{x}^{2lr}+\int_t^\tau \vert W_A(s)\vert^{2lr} ds+\int_t^{\tau}\int_t^u\abs{W_A(s)}^{2lr} ds du\right).
\]
Finally, we obtain
\begin{eqnarray*}
\E\left[ e^{p\int_0^T\vert X_s^{t,x}\vert^{2lr}ds}\right]
 &\leq& Ce^{C\abs{x}^{2lr}}\E\left[ e^{C\int_t^T\left(\vert \Gamma_s^{t,x}\vert^{2lr}
 +\vert W_A(s)\vert^{2lr}\right)\,ds}\right]\\
&\leq& Ce^{C\abs{x}^{2lr}}\E\left[ e^{C\int_t^T\vert W_A(s)\vert^{2lr}\,ds+C\int_t^T\int_t^{\tau} \vert W_A(s)\vert^{2lr}\,dsd\tau+C\int_t^T\int_t^{\tau}\int_t^u \vert W_A(s)\vert^{2lr}\,dsdud\tau}\right]\\
&\leq& Ce^{C\abs{x}^{2lr}}\E\left[ e^{C\int_t^T\vert W_A(s)\vert^{2lr}\,ds}\right]^{1/3}\E\left[ e^{C\int_t^T\int_t^{\tau} \vert W_A(s)\vert^{2lr}\,dsd\tau}\right]^{1/3}\\
& & \times \E\left[e^{C\int_t^T\int_t^{\tau}\int_t^u \vert W_A(s)\vert^{2lr}\,dsdud\tau}\right]^{1/3}.
\end{eqnarray*}
By using Jensen inequality as in the proof of Proposition \ref{prop-exist-bsde}, we can show that
\begin{equation*}
 \E\left[ e^{C\int_t^T\vert W_A(s)\vert^{2lr}\,ds}\right] \leq \frac{C}{T-t}\int_t^T\E\left[ e^{C\vert W_A(s)\vert^{2lr}}\right]\,ds\leq C
\end{equation*}
because $W_A(s)$ is a centered Gaussian random variable with a bounded covariance operator. By same arguments, we are able to show that
\begin{equation*}
 \E\left[ e^{C\int_t^T\int_t^{\tau} \vert W_A(s)\vert^{2lr}\,dsd\tau}\right]\leq C \quad \textrm{and} \quad \E\left[e^{C\int_t^T\int_t^{\tau}\int_t^u \vert W_A(s)\vert^{2lr}\,dsdud\tau}\right] \leq C,
\end{equation*}
which achieved the proof of (\ref{exponential moment sup X2}).

The end of the proof goes on like the proof of Proposition \ref{prop-exist-bsde} by assuming the following lemma
\ref{lemma-iterate-pert}, analogous of Lemma \ref{lemma-iterate}.
\qed

\begin{lemma}\label{lemma-iterate-pert}
 Under the assumptions of Proposition \ref{prop-exist-bsde}, we have
\begin{equation*}
\vert Z_\tau^M\vert\leq A_n+B_n\vert X_\tau\vert^r,
\end{equation*}
with $(A_n,B_n)_{n\in\N}$ defined by recursion: $B_0=0$, $A_0$ defined before,
\begin{eqnarray*}
 A_{n+1} &=& C(1+A_n^{lr}),\\
 B_{n+1} &=& C,
\end{eqnarray*}

where $C$ is a constant that does not depend on $M$.
\end{lemma}
\noindent {\bf Proof.} The only difference with the proof of Lemma
\ref{lemma-iterate} arise when estimating
 $\E_\tau^{\Q^M}\vert X_s^{t,x}\vert^r$, for $s\in[\tau,T]$. We have only to notice
that this time in $(\Omega,\calf,\Q^M)$, for $s\in[\tau,T]$,
$X_s^{t,x}$ solves the following equation in mild form
\begin{align*}
 X_s^{t,x}&=e^{(s-\tau)A}X_\tau^{t,x}+\int_\tau^s e^{(s-r)A}B\,d\tilde{W_r}\\
&+\int_\tau^s e^{(s-r)A}B\left(F(r,X_r^{t,x})+ 
\nabla_z\psi_M(r,X_r^{t,x},Y^{M,t,x}_r,Z^{M,t,x}_r)\right)\,dr.
\end{align*}
So arguing as in the proof of Lemma \ref{lemma-iterate} we get the conclusion.
\qed

We only mention that, as in \cite{Ri1}, the results
contained in Propositions \ref{prop-exist-bsde} and \ref{prop-exist-bsde-pert}
could be stated
under slightly weaker assumptions than Hypothesis 3.1: we could threat the case $rl=1$ when $T$ is small enough. Nevertheless any applications on HJB equations follow
by weakening the assumptions in that direction.

\section{Differentiability with respect to the initial datum in the FBSDE}

\label{sec-diffle}
In this section we consider regular dependence on the initial datum of the
perturbed Ornstein-Uhlenbeck process $X$ for the Markovian BSDE (\ref{bsde}), namely we consider once again the following decoupled forward backward system
\begin{equation}
\left\{
\begin{array}
[c]{l}%
dX_\tau^{t,x}  =AX_\tau^{t,x} d\tau+F(\tau,X_\tau^{t,x}) d\tau+BdW_\tau
,\text{ \ \ \ }\tau\in\left[  t,T\right], \\ \dis
X_{\tau}^{t,x} =x, \text{ \ \ \ }\tau\in\left[  0,t\right],\\ \dis
 dY_\tau^{t,x}=-\psi(\tau, X^{t,x}_\tau,Y_\tau^{t,x},Z_\tau^{t,x})\;d\tau+Z^{t,x}_\tau\;dW_\tau,
 \qquad \tau\in [0,T],
  \\\dis
  Y_T^{t,x}=\phi(X_T^{t,x}).
\end{array}
\right.  \label{fbsde}%
\end{equation}
Beside Hypotheses \ref{ip su AB} and \ref{ip-psi-phi} on the coefficients, we assume
the following hypothesis.
\begin{hypothesis}\label{ip-psi-phi-diffle}
$ $
\begin{enumerate}
\item For every $\tau\in[0,T]$, the map
$(x,y,z)\mapsto\psi(\tau,x,y,z)$ belongs to $\calg^{1,1,1}(H\times\R\times\Xi^{\ast},\R)$, and by
Hypothesis \ref{ip-psi-phi}, 
\[
\abs{\nabla_x \psi(\tau,x,y,z)}\leq \left(C+\beta\vert x\vert^r\right),
\quad \abs{\nabla_y \psi(\tau,x,y,z)}\leq K_{\psi_y},
\quad \abs{\nabla_z \psi(\tau,x,y,z)}\leq \left(C+\gamma\vert z\vert^l\right),
\]
$\forall \tau\in[0,T],\,
\forall x\in H,\,\forall y\in\R,\,\forall z\in\Xi^{\ast}$.
\item $\phi\in\calg^1(H,\R)$ and by Hypothesis \ref{ip-psi-phi},
$$\abs{\nabla_x \phi(x)} \leq \left(C+\alpha \abs{x}^r\right), \qquad \forall x \in H.$$
\end{enumerate}
\end{hypothesis}
The following result is proved by Fuhrman and Tessitore in \cite{fute}.

\begin{proposition}
\label{theo-Xdiffle}
 Assume Hypothesis \ref{ip su AB} holds true. Then the map $(t,x) \mapsto X^{t,x}$ belongs to $\calg^{0,1}([0,T]
\times H;\cals^p)$ for all $p>1$. Moreover, we have, for every $x,h \in H$,
$$\Vert \nabla_x X_{\tau}^{t,x} h \Vert_{\mathcal{S}^p} \leq C_p \vert h \vert.$$
\end{proposition}

We are now able to give the main result of this section.

 \begin{theorem}\label{teo-diffle}
Assume Hypotheses \ref{ip su AB}, \ref{ip-psi-phi} and \ref{ip-psi-phi-diffle} hold true.
Then the map $(t,x)\mapsto (Y^{t,x},Z^{t,x})$ belongs to $\calg^{0,1}([0,T]
\times H;\cals^p\times\calm^p)$ for each $p>1$. Moreover for every $x,h\in H$
the directional derivative process $(\nabla_x Y_\tau^{t,x},\nabla_x Z^{t,x}_\tau)_{\tau\in[0,T]}$
solves the following BSDE: for $\tau\in[0,T],$
\begin{align}
\nonumber \nabla_x Y_\tau^{t,x}h=&\nabla_x\phi(X_T^{t,x})\nabla_xX_T^{t,x}h
+\int_\tau^T\nabla_x\psi(s, X^{t,x}_s,Y_s^{t,x},Z_s^{t,x})\nabla_xX_s^{t,x}h\;ds\\ \nonumber
&+\int_\tau^T\nabla_y\psi(s, X^{t,x}_s,Y_s^{t,x},Z_s^{t,x})\nabla_xY_s^{t,x}h\;ds
+\int_\tau^T\nabla_z\psi(s, X^{t,x}_s,Y_s^{t,x},Z_s^{t,x})\nabla_xZ_s^{t,x}h\;ds\\ 
&-\int_\tau^T \nabla_x Z^{t,x}_sh\;dW_s  \label{BSDE diff}
\end{align}
and there exist two constants $C$ and $C_p$ such that 
$$\abs{\nabla_x Y_{\tau}^{t,x} h } \leq C\left(1+\abs{X_{\tau}^{t,x}}^r\right)\abs{h},\quad \quad\Vert \nabla_x Z^{t,x} h \Vert_{\mathcal{M}^p} \leq C_p e^{C_p \abs{x}^{2rl}}\vert h \vert.$$
\end{theorem}
\noindent {\bf Proof.}
Firstly, we will show the continuity of the map $(t,x) \mapsto (Y^{t,x},Z^{t,x})$. We fix $(t,x) \in [0,T] \times H$ and we consider $(t',x') \in [0,T] \times H$ such that $t' \rightarrow t$ and $x' \rightarrow x$. Let us denote
$$\delta Y := Y^{t,x}-Y^{t',x'} \quad \textrm{and} \quad \delta Z :=Z^{t,x}-Z^{t',x'}.$$
The usual linearization trick gives us that $(\delta Y,\delta Z)$ is the solution of the BSDE
\begin{align*}
 \delta Y_s=&\phi(X_T^{t,x})-\phi(X_T^{t',x'})-\int_s^T \delta Z_u dW_u\\
& + \int_s^T \left[\psi(u,X_u^{t,x},Y_u^{t,x},Z_u^{t,x})-\psi(u,X_u^{t',x'},Y_u^{t,x},Z_u^{t,x})+U_u \delta Y_u+\langle V_u, \delta Z_u\rangle_{\Xi^{\ast}}\right]\, du,
\end{align*}
with
\[
U_u=\left\lbrace
\begin{array}{ll}
\dfrac{\psi(u,X_u^{t',x'}, Y^{t,x}_u,Z^{t,x}_u)
-\psi(u,X_u^{t',x'}, Y^{t',x'}_u,Z^{t,x}_u) }{Y_u^{t,x}-Y_u^{t',x'}} &\text{if }Y_u^{t,x}-Y_u^{t',x'}\neq 0\\
0&\text{if }Y_u^{t,x}-Y_u^{t',x'}= 0
\end{array}
 \right.\]
and\[
 V_u=\left\lbrace
 \begin{array}{ll}
 \dfrac{\psi(u,X_u^{t',x'}, Y^{t',x'}_u,Z^{t,x}_u)
 -\psi(u,X_u^{t',x'}, Y^{t',x'}_u,Z^{t',x'}_u) }{\abs{Z^{t,x}_u-Z^{t',x'}_u}^2}\left(Z^{t,x}_u-Z^{t',x'}_u\right) &\text{if }Z^{t,x}_u-Z^{t',x'}_u\neq 0\\
 0&\text{if }Z^{t,x}_u-Z^{t',x'}_u= 0.
 \end{array}
 \right.
   \]
Thanks to Hypothesis \ref{ip-psi-phi}, we remark that $\vert U_u \vert \leq K_{\psi_y}$ and 
$$\vert V_u \vert \leq C(1+\vert Z_u^{t,x} \vert^l + \vert Z_u^{t',x'} \vert^l) \leq C(1+\vert X_u^{t,x} \vert^{rl} + \vert X_u^{t',x'} \vert^{rl}).$$
A mere extension of Proposition 3.6 in \cite{WaRaCh} gives us a stability result: for all $p>1$
\begin{eqnarray*}
& &\Vert \delta Y \Vert_{\mathcal{S}^p}^p+\Vert \delta Z \Vert_{\mathcal{M}^p}^p\\
& \leq &C_p \mathbb{E} \left[ e^{4p \int_0^T \vert V_s \vert^2ds} \vert \phi(X_T^{t,x})-\phi(X_T^{t',x'}) \vert^{2p}\right]\\
&&+C_p \mathbb{E} \left[\left(\int_0^T e^{4 \int_0^s \vert V_u \vert^2du}\vert \psi(s,X_s^{t,x},Y_s^{t,x},Z_s^{t,x})-\psi(s,X_s^{t',x'},Y_s^{t',x'},Z_s^{t',x'}) \vert ds\right)^p\right]\\
& \leq &C_p \mathbb{E} \left[ e^{C_p \int_0^T \left(\vert X_s^{t,x} \vert^{2rl}+\vert X_s^{t',x'} \vert^{2rl}\right)ds}(1+\vert X_T^{t,x} \vert^{2pr}+\vert X_T^{t',x'} \vert^{2pr}) \vert  X_T^{t,x}-X_T^{t',x'}\vert^{2p}\right]\\
&&+C_p \int_0^T\mathbb{E} \left[ e^{C_p \int_0^s \left(\vert X_u^{t,x} \vert^{2rl}+\vert X_u^{t',x'} \vert^{2rl}\right)du} (1+\vert X_s^{t,x} \vert^{2pr}+\vert X_s^{t',x'} \vert^{2pr}) \vert  X_s^{t,x}-X_s^{t',x'}\vert^{2p} \right]ds.
\end{eqnarray*}
By using H\"older theorem, Proposition \ref{esistenza sde}, estimate (\ref{exponential moment sup X2}) and classical stability results for the solution of the forward equation, we show that the right term in the last inequality tends to $0$ when $t' \rightarrow t$ and $x' \rightarrow x$. So we have that $(t,x) \mapsto (Y^{t,x},Z^{t,x})$ is continuous in $\mathcal{S}^p \times \mathcal{M}^p$ for all $p>1$.

For the differentiability, we will follow the proof of Proposition 12 in \cite{BriFu}. Firstly, let us remark that, thanks to Hypothesis \ref{ip-psi-phi-diffle} and Proposition \ref{prop-exist-bsde-pert}, 
$$\vert \nabla_z \psi (s,X_s^{t,x},Y_s^{t,x},Z_s^{t,x}) \vert \leq C(1+\vert Z_s^{t,x} \vert^l) \leq C(1+\vert X_s^{t,x} \vert^{rl}),$$
and, thanks to Propositions \ref{esistenza sde} and \ref{theo-Xdiffle}, estimate (\ref{exponential moment sup X2}) and Hypothesis \ref{ip-psi-phi-diffle}, for all $p>1$, for all $c>0$ and for all $h \in H$,
$$ \mathbb{E} \left[ e^{c\int_0^T \vert X_s^{t,x} \vert^{2rl}ds}\left(\vert \nabla_x \phi(X_T^{t,x})\nabla_x X_T^{t,x} h \vert^p+ \int_0^T \vert \nabla_x \psi(s,X_s^{t,x},Y_s^{t,x},Z_s^{t,x})\nabla_x X_s^{t,x} h \vert^p ds \right)\right]<+\infty.$$
So, it follows from a mere generalization of Theorem 4.1 in \cite{WaRaCh} that BSDE (\ref{BSDE diff}) has a unique solution which belongs to $\mathcal{S}^p \times \mathcal{M}^p$ for all $p>1$. Now, let us fix $(t,x) \in [0,T] \times H$. We remove parameters $t$ and $x$ for notational simplicity. For $\varepsilon>0$, we set $X^{\varepsilon}:=X^{t,x+\varepsilon h}$, where $h$ is some vector in $H$, and we consider $(Y^{\varepsilon},Z^{\varepsilon})$ the solution in $\mathcal{S}^p \times \mathcal{M}^p$ to the BSDE
$$Y_t^{\varepsilon} = \phi(X_T^{\varepsilon})+\int_t^T \psi(s,X_s^{\varepsilon},Y_s^{\varepsilon},Z_s^{\varepsilon})ds-\int_t^T Z_s^{\varepsilon}dW_s.$$
When $\varepsilon \rightarrow 0$, $(X^{\varepsilon},Y^{\varepsilon},Z^{\varepsilon}) \rightarrow (X,Y,Z)$ in  $\mathcal{S}^p \times \mathcal{S}^p \times \mathcal{M}^p$, for all $p>1$. We also denote $(G,N)$ the solution to the BSDE (\ref{BSDE diff}). We have to prove that the directional derivative of the map $(t,x) \mapsto (Y^{t,x}, Z^{t,x})$ in the direction $h \in H$ is given by $(G,N)$. Let us consider $U^{\varepsilon}:=\varepsilon^{-1}(Y^{\varepsilon}-Y)-G$, $V^{\varepsilon}:=\varepsilon^{-1}(Z^{\varepsilon}-Z)-N$. We have
\begin{eqnarray*}
 U^{\varepsilon}_t &=& \frac{1}{\varepsilon} \left(\phi(X_T^{\varepsilon})-\phi(X_T)\right)-\nabla_x \phi(X_T)\nabla_x X_Th\\
 && +\frac{1}{\varepsilon} \int_t^T \left(\psi(s,X_s^{\varepsilon},Y_s^{\varepsilon},Z_s^{\varepsilon})-\psi(s,X_s,Y_s,Z_s)\right)ds-\int_t^T V_s^{\varepsilon} dW_s\\
 && -\int_t^T \nabla_x \psi(s,X_s,Y_s,Z_s) \nabla_x X_s h ds - \int_t^T \nabla_y \psi(s,X_s,Y_s,Z_s) G_s ds\\
 && - \int_t^T \nabla_z \psi(s,X_s,Y_s,Z_s) N_s ds.
\end{eqnarray*}
As in the proof of Proposition 12 in \cite{BriFu}, we use the fact that $\psi(s,\cdot,\cdot,\cdot)$ belongs to $\mathcal{G}^{1,1,1}$ and so we can write
\begin{eqnarray*}
 &&\frac{1}{\varepsilon} \left(\psi(s,X_s^{\varepsilon},Y_s^{\varepsilon},Z_s^{\varepsilon})-\psi(s,X_s,Y_s,Z_s)\right)\\
 &=& \frac{1}{\varepsilon} \left(\psi(s,X_s^{\varepsilon},Y_s,Z_s)-\psi(s,X_s,Y_s,Z_s)\right)+A_s^{\varepsilon} \frac{Y_s^{\varepsilon}-Y_s}{\varepsilon}+B_s^{\varepsilon}\frac{Z_s^{\varepsilon}-Z_s}{\varepsilon},
\end{eqnarray*}
where $A_s^{\varepsilon} \in L(\mathbb{R}, \mathbb{R})$ and $B_s^{\varepsilon} \in L(\Xi^{\ast},\mathbb{R})$ are defined by
\begin{eqnarray*}
 A^{\varepsilon}_s y &:=& \int_0^1 \nabla_y \psi(s,X_s^{\varepsilon},Y_s+\alpha(Y_s^{\varepsilon}-Y_s),Z_s)yd\alpha, \quad \forall y \in \mathbb{R},\\
 B^{\varepsilon}_s z &:=& \int_0^1 \nabla_z \psi(s,X_s^{\varepsilon},Y_s^{\varepsilon},Z_s+\alpha(Z_s^{\varepsilon}-Z_s))zd\alpha, \quad \forall z \in \Xi^{\ast}.
\end{eqnarray*}
Then $(U^{\varepsilon},V^{\varepsilon})$ solves the BSDE:
$$U^{\varepsilon}_t=\zeta^{\varepsilon}+ \int_t^T (A^{\varepsilon}_s U^{\varepsilon}_s+B^{\varepsilon}_s V^{\varepsilon}_s)ds+\int_t^T (P^{\varepsilon}_s+Q^{\varepsilon}_s+R^{\varepsilon}_s)ds-\int_t^T V^{\varepsilon}_s dW_s,$$
where we have set
\begin{eqnarray*}
 \zeta^{\varepsilon} &:=& \varepsilon^{-1}(\phi(X_T^{\varepsilon})-\phi(X_T))-\nabla_x \phi(X_T)\nabla_x X_Th,\\
 P^{\varepsilon}_s &:=& (A_s^{\varepsilon}-\nabla_y \psi(s,X_s,Y_s,Z_s))G_s,\\
 Q^{\varepsilon}_s &:=& (B_s^{\varepsilon}-\nabla_z \psi(s,X_s,Y_s,Z_s))N_s,\\
 R_s^{\varepsilon} &:=& \varepsilon^{-1}( \psi(s,X_s^{\varepsilon},Y_s,Z_s)- \psi(s,X_s,Y_s,Z_s))-\nabla_x \psi(s,X_s,Y_s,Z_s)\nabla_x X_sh.
\end{eqnarray*}
It follows from Hypothesis \ref{ip-psi-phi-diffle} and estimate (\ref{estimate-Z-r-Fneqzero}) that
\begin{eqnarray*}
 && \vert A^{\varepsilon}_s\vert \leq C, \quad \vert B^{\varepsilon}_s \vert \leq C(1+\vert X_s \vert^{rl}+\vert X_s^{\varepsilon} \vert^{rl}), \\
 && \vert P^{\varepsilon}_s \vert \leq C \vert G_s \vert, \quad \vert Q^{\varepsilon}_s \vert \leq C(1+ \vert X_s \vert^{rl}+\vert X_s^{\varepsilon} \vert^{rl})\vert N_s \vert.
\end{eqnarray*}
We have, once again from a mere generalization of Proposition 3.6 in \cite{WaRaCh}, 
\begin{eqnarray*}
\Vert U^{\varepsilon} \Vert_{\mathcal{S}^p}+ \Vert V^{\varepsilon} \Vert_{\mathcal{M}^p} &\leq & C_p \mathbb{E} \left[ e^{C_p \int_0^T \vert X_s \vert^{rl}+\vert X_s^{\varepsilon} \vert^{rl} ds} \left( \vert \zeta^{\varepsilon} \vert^{2p} + \int_0^T \vert P^{\varepsilon}_s\vert^{2p}+\vert Q^{\varepsilon}_s\vert^{2p}+\vert R^{\varepsilon}_s\vert^{2p} ds\right) \right].
\end{eqnarray*}
By using H\"older inequality and the estimate (\ref{estimate-Z-r-Fneqzero}), previous inequality becomes
\begin{eqnarray*}
\Vert U^{\varepsilon} \Vert_{\mathcal{S}^p}+ \Vert V^{\varepsilon} \Vert_{\mathcal{M}^p} &\leq & C_p e^{C_p \left(\abs{x}^{rl}+\abs{x+\varepsilon h}^{rl}\right)} \mathbb{E} \left[ \vert \zeta^{\varepsilon} \vert^{4p} + \int_0^T \vert P^{\varepsilon}_s\vert^{4p}+\vert Q^{\varepsilon}_s\vert^{4p}+\vert R^{\varepsilon}_s\vert^{4p} ds \right].
\end{eqnarray*}
By using a uniform integrability argument, the right hand side of the previous inequality tends to $0$ as $\varepsilon \rightarrow 0$ in view of the regularity and the growth of $\phi$ and $\psi$.

The proof that maps $x \mapsto (\nabla_x Y^{t,x}h,\nabla_xZ^{t,x}h)$ and $h \mapsto (\nabla_x Y^{t,x}h,\nabla_xZ^{t,x}h)$ are continuous (for every $h$ and $x$ respectively) comes once again from a mere generalization of Proposition 3.6 in \cite{WaRaCh}.

To finish the proof, it remains to prove the growth estimate on $\abs{\nabla_x Y^{t,x} h}$ and $\Vert \nabla_x Z^{t,x} h \Vert_{\mathcal{M}^p}$. Let us begin with the first one. Thanks to the estimate on $Z^{t,x}$ given by proposition \ref{prop-exist-bsde-pert}, we have
$$\abs{\nabla_z \psi(s,X_s^{t,x},Y_s^{t,x},Z_s^{t,x})} \leq C(1+\abs{X_s^{t,x}}^{rl}).$$
Now the result (\ref{exponential moment sup X2}) shows us that Novikov's condition is fulfilled and so we are able to use Girsanov's theorem in (\ref{BSDE diff}): there exists a probability $\mathbb{Q}$, equivalent to the original one $\mathbb{P}$, such that $\tilde{W}_{\tau}:=W_{\tau}-\int_0^{\tau} \nabla_z \psi(s,X_s^{t,x},Y_s^{t,x},Z_s^{t,x})ds$ is a Wiener process under $\mathbb{Q}$. We obtain 
\begin{align*}
 \nabla Y^{t,x}_\tau h&=\E_\tau^{\Q}\left[ e^{\int_\tau^T
\nabla_y\psi(u,X_u^{t,x},Y^{t,x}_u,Z^{t,x}_u)du}\nabla\phi(X_T^{t,x})\nabla X_T^{t,x}h\right.\\ \nonumber
&\left.+\int_\tau^Te^{\int_\tau^T
\nabla_y\psi(u,X_u^{t,x},Y^{t,x}_u,Z^{t,x}_u)du}\
\nabla_x\psi(s,X_s^{t,x},Y^{t,x}_s,Z^{t,x}_s)\nabla X_s^{t,x}hds\right],
\end{align*}
and, by using assumptions \ref{ip-psi-phi-diffle} and the fact that $\nabla X^{t,x}$ is bounded,
$$\abs{\nabla Y^{t,x}_\tau h} \leq C\left(1+ \E_{\tau}^{\Q} \left[ \abs{X_T^{t,x}}^r \right] +\int_{\tau}^T \E_{\tau}^{\Q} \left[ \abs{X_s^{t,x}}^r \right]ds \right)\abs{h}.$$
Then, arguing as in the proof of Lemma \ref{lemma-iterate-pert}, we obtain that $\E_{\tau}^{\Q} \left[  \abs{X_s^{t,x}}^r \right] \leq C\left(1+\abs{X_{\tau}^{t,x}}^r\right)$ and, finally,
$$\abs{\nabla Y^{t,x}_\tau h} \leq C\left(1+\abs{X_{\tau}^{t,x}}^r \right)\abs{h}.$$
For the estimate on  $\Vert \nabla_x Z^{t,x} h \Vert_{\mathcal{M}^p}$, we just have to use a mere generalization of Proposition 3.6 in \cite{WaRaCh}.
\qed

\section{Probabilistic solution of a semilinear PDE in infinite dimension: the differentiable data case}
\label{sezionePDEdiffle}
The aim of this section is to present existence and uniqueness results for the solution of
a semilinear Kolmogorov equation with the nonlinear term which is superquadratic 
with respect to the $B$-derivative and with final datum not necessarily bounded, in 
the case of differentiable coefficients.

More precisely, let $\mathcal L_t $ be the generator of the transition
semigroup $(P_{t,\tau})_{\tau\in[t,T]}$, that is, at least formally, 
$$
(\call_t f)(x)=\frac{1}{2}(Tr BB^* \nabla^2 f)(x)+\<Ax,\nabla f(x)\>+\<F(t,x),\nabla f(x)\>.
$$
Let us consider the following equation
\begin{equation}
\left\{
\begin{array}
[c]{l}%
\frac{\partial v}{\partial t}(t,x)=-\mathcal{L}_tv\left(  t,x\right)
+\psi\left(t,x,v(t,x),  \nabla^{B}v\left(  t,x\right)  \right)  ,\text{ \ \ \ \ }t\in\left[  0,T\right]
,\text{ }x\in H,\\
v(T,x)=\phi\left(  x\right).
\end{array}
\right.  \label{Kolmo}%
\end{equation}
In the following we introduce the notion of mild solution for the non linear Kolmogorov
equation (\ref{Kolmo}) (see also \cite{DP1} and \cite{fute}, or \cite{Mas} for the definition of
mild solution when $\psi$ depends only on $\nabla^{B}v$ and not on $\nabla v$).

\noindent Notice that, by Proposition \ref{esistenza sde}, if $\phi$ satisfies
Hypothesis \ref{ip-psi-phi}, point 3, or more generally if $\phi$ is a continuous
function with polynomial growth,
by $L^p(\Omega,C([0,T]))$-integrability of any order $p$ of the Markov process $X^{t,x}$,
given by Proposition \ref{esistenza sde}, we have that
\[
 P_{t,\tau}[\phi](x)=\E \left[\phi(X^{t,x}_\tau)\right]
\]
 is well defined. Since $\mathcal{L}_t$ is (formally) the generator of
$(P_{t,\tau})_{\tau\in[t,T]}$, the variation of constants formula for  equation (\ref{Kolmo}) gives us:%
\begin{equation}
v(t,x)=P_{t,T}\left[  \phi\right]  \left(  x\right)  +\int_{t}^{T}%
P_{t,s}\left[  \psi\left(s,\cdot,v\left(  s,\cdot\right),\nabla^B 
v\left(  s,\cdot\right)\right)  \right] (  x)  ds,
\text{\ \ }t\in\left[  0,T\right]  ,\text{ }x\in H. \label{solmildkolmo}%
\end{equation}
We will use this formula to define the notion of mild
solution for the non linear Kolmogorov equation (\ref{Kolmo}); before giving the definition we
have also to
introduce some spaces of continuous functions, where we will look for the
solution of (\ref{Kolmo}).

\noindent We consider the space $C_{b}^{s}\left(
H,\Xi^{\ast}\right)  $ of mappings $L:H\rightarrow\Xi^{\ast}$ such that
for every $\xi\in\Xi$, $L\left(  \cdot\right)  \xi\in C_{b}\left(  H\right)
$, where $C_{b}\left(  H\right)$ denotes the space of bounded continuous functions from $H$ to $\mathbb{R}$. The space $C_{b}^{s}\left(  H,\Xi^{\ast}\right)  $ turns out to be a Banach
space if it is endowed with the norm
\begin{equation*}
\left\|  L\right\|  _{C_{b}^{s}\left(H,  \Xi^{\ast}\right)  }=\sup_{x\in
H}\left|  L\left(  x\right)  \right|  _{\Xi^{\ast}}. 
\end{equation*}
Besides $C_{b}^{s}\left(  H,\Xi^{\ast}\right)  $ we consider also the linear
space $C_{k}^{s}\left(  H,\Xi^{\ast}\right)  $ of mappings $L:H\rightarrow
\Xi^{\ast}$ such that for every $\xi\in\Xi$, $L\left(  \cdot\right)  \xi\in
C_{k}\left(  H\right)$, where $C_{k}\left(  H\right)$ denotes the space of continuous functions from $H$ to $\mathbb{R}$ with a polynomial growth of degree $k$. The linear space $C_{k}^{s}\left(  H,\Xi^{\ast
}\right)  $ turns out to be a Banach space if it is endowed with the norm%

\[
\left\|  L\right\|  _{C_{k}^{s}\left(  H,\Xi^{\ast}\right)  }=\sup_{x\in
H}\frac{\left|  L\left(  x\right)  \right|  _{\Xi^{\ast}}}{\left(  1+\left|
x\right|  _{H}^{2}\right)  ^{k/2}}.
\]

We are now able to give the definition of a mild solution of \eqref{Kolmo}.
\begin{definition}
\label{defsolmildkolmo}Let $r\geq 0$.
We say that a function $v:\left[  0,T\right]  \times H\rightarrow\mathbb{R}$ is a mild
solution of the non linear Kolmogorov equation (\ref{Kolmo}) if the following
are satisfied:

\begin{enumerate}
\item $v\in C_{r+1}\left(  \left[  0,T\right]  \times H\right)  $;

\item $\nabla^{B}  v\in C^{s}_{r}\left(  \left[
0,T\right)  \times H,\Xi^{\ast}\right)  $,
in particular this means that for
every $t\in\left[  0,T\right)  $, $v\left(  t,\cdot\right)$ is $B$-differentiable
and the derivative has polynomial growth of order $r$;

\item  equality (\ref{solmildkolmo}) holds.
\end{enumerate}
\end{definition}

Notice that the differentiability required at point 2 is the minimal request in order to make equality
(\ref{solmildkolmo}) work. In the case of differentiable data $\psi$ and $\phi$, in addition to
differentiability of the nonlinear term $F$ in the forward equation (\ref{ornstein-pert}),
we look for a solution
$v$ differentiable with respect to $x$ in all directions. In this case
$\nabla^B v=\nabla v\,B$ and saying that a function $v:[0,T]\times H\rightarrow\R$
admits a G\^ateaux derivative $\nabla v\in C_k([0,T]\times H,H^{\ast})$
is equivalent to ask $v\in \calg^{0,1}([0,T]\times H)$ such that the operator norm
of $\nabla v(t,x)$ has polynomial growth of order $k$ with respect to $x$.
So, in this part we will prove the existence of a mild solution according to the following stronger definition:
\begin{definition}
\label{defsolmildkolmosmooth}Let $r\geq 0$. We say that a
function $v:\left[  0,T\right]  \times H\rightarrow\mathbb{R}$ is a mild
solution of the non linear Kolmogorov equation (\ref{Kolmo}) if the following
are satisfied:

\begin{enumerate}
\item $v\in C_{r+1}\left(  \left[  0,T\right]  \times H\right)  $;

\item for every $t\in[0,T]$, $v(t,\cdot)$ is differentiable in $H$ and the derivative
has polynomial growth with respect to $x$, more precisely
$v\in\calg^{0,1}([0,T]\times H)$ and $\forall\,h\in H$
\[\sup
_{t\in\left[  0,T\right]  }\sup_{x\in H} \dfrac{\left\vert \nabla_x v\left(  t,x\right)h
\right\vert }{(1+\vert x\vert^2)^{r/2}}  <\infty;
\]
\item  equality (\ref{solmildkolmo}) holds.
\end{enumerate}
\end{definition}
We notice that we will take in the following the same index $r$ than in Hypothesis \ref{ip-psi-phi},
so this index is related to the growth of $\phi$ and $\psi$ with respect to $x$.

Existence and uniqueness of a mild solution of equation (\ref{Kolmo})
is related to the study of the
forward-backward system given by the perturbed Ornstein-Uhlenbeck
process $X^{t,x}$ defined in (\ref{ornstein-pert}) and by the BSDE (\ref{bsde}). 
We will show that, if we define 
\begin{equation*}
v(t,x):=Y_t^{t,x},
\end{equation*}
with $(Y^{t,x},Z^{t,x})$ the solution of the BSDE (\ref{bsde}), then it turns out that $v$ is the unique mild solution of equation (\ref{Kolmo}),
and $\nabla^B v(t,x)=Z_t^{t,x}$. On the coefficients $\psi$, $\phi$ and $F$ of equation (\ref{Kolmo}), which are the same appearing in the backward equation in the system (\ref{fbsde}) and on the non linear term
of the forward equation in the system (\ref{fbsde}), we make
differentiability assumptions contained in Hypothesis \ref{ip-psi-phi-diffle}.

Notice that we are working with a function $\psi$ that can have a quadratic ($l=1$)
or a superquadratic growth ($l>1$) with respect to $z$.
Moreover, $\psi$ and $\phi$ are unbounded and can have some polynomial growth with respect to $x$,
though this growth is forced to decrease as the growth with respect to $z$ increases,
see again Hypothesis \ref{ip-psi-phi}.
So the result we are going to obtain improves Theorem 15 in \cite{BriFu},
where it is considered the quadratic case for $\psi$ with respect to $z$
and a bounded final datum,
and also Theorem 4.1 in \cite{Mas3}, where the superquadratic case is considered
in the case of a bounded final datum together with some smoothing properties for the transition
semigroup of the forward equation. Notice that we will require similar
smoothing properties in the next section, when we will
remove differentiability assumptions on the coefficients.

\begin{theorem}\label{teokolmosmooth}
Assume that Hypotheses \ref{ip su AB}, \ref{ip-psi-phi}, \ref{ip-psi-phi-diffle}
hold true. Then, according to definition \ref{defsolmildkolmosmooth},
equation (\ref{Kolmo}) admits a unique mild solution. This solution satisfies
\[
\vert v(t,x)\vert\leq C(1+\vert x\vert^{r+1}),\qquad
\vert \nabla^B v(t,x)\vert\leq C(1+\vert x\vert^{r}).
\qquad 
\]
\end{theorem}

\noindent {\bf Proof.} The proof is substantially based on estimate (\ref{estimate-Z-r})
and on section \ref{sec-diffle} where differentiability of the FBSDE (\ref{fbsde})
in the case of differentiable coefficients is investigated.
Since we assume that coefficients are differentiable, by Theorem \ref{teo-diffle} 
$Y^{t,x}$ is differentiable with respect to $x$. We set $v(t,x):=Y^{t,x}_t$: 
notice that as usual $Y^{t,x}_t$ is deterministic.
As in Lemma 6.3 in \cite{fute},
we can prove that $\forall \xi \in \Xi$ and $\forall s\in[t,T]$ the joint quadratic variation
\[
 \langle v(s,X_s^{t,x}),\int_t^s<\xi, dW_{\tau}>\rangle
=\int_t^s \nabla_x v(\tau,X_\tau^{t,x})B\xi\,d\tau.
\]
Since $v(s,X_s^{t,x})=Y_s^{t,x}$, from the BSDE in (\ref{fbsde}) we get that
$\forall \xi \in \Xi$ and $\forall s\in[t,T]$ the joint quadratic variation is equal to
\[
 \langle v(s,X_s^{t,x}),\int_t^s<\xi, dW_{\tau}>\rangle=\int_t^s Z_\tau^{t,x}\xi\,d\tau.
\]
This gives the identification, for a.a. $\tau\in[t,T]$,
\begin{equation}\label{identif-Z}
 Z_\tau^{t,x}=\nabla_x v(\tau,X_\tau^{t,x})B \quad \P\textrm{-a.s. }.
\end{equation}
With this identification in hand, the proof goes on in a quite standard way: see e.g.
the proof of Theorem 6.2 in the pioneering paper \cite{fute} for the study of BSDEs and related PDEs in infinite dimension. We give here a sketch of the proof for the reader convenience.

\emph{Existence.} Let us consider $(Y^{t,x},Z^{t,x})$ the solution of the BSDE (\ref{fbsde}),
which in integral form is given by
\begin{equation*}
  Y_s^{t,x}+\int_s^T Z^{t,x}_\tau\;dW_\tau=\phi(X_T^{t,x})+\int_s^T\psi(\tau, X^{t,x}_\tau,Y_\tau^{t,x},Z_\tau^{t,x})\;d\tau
\end{equation*}
Taking expectation, setting $s=t$ and using (\ref{identif-Z}) we get 
the existence of a mild solution according to definition
\ref{defsolmildkolmosmooth}: notice that the growth of
$\nabla^B v$ comes from estimates on $Z$ in Propositions \ref{prop-exist-bsde}
and \ref{prop-exist-bsde-pert}, namely see estimates (\ref{estimate-Z-r}) and (\ref{estimate-Z-r-Fneqzero}).
For what concerns the estimate on $v$, we can mimic the proof of Proposition 2.5 in \cite{MR}, and then obtain the
desired polynomial growth for $v$ with respect to $x$:
\begin{equation*}
  \vert v(t,x)\vert:=\vert Y^{t,x}_t\vert\leq C\left(1+\vert x\vert^{r+1}\right).
 \end{equation*}

\medskip

\emph{Uniqueness.} Let $u$ be a mild solution of equation 
(\ref{Kolmo}): by the Markov property of the process
$X^{t,x}$, we have, $\forall\,s\in[t,T]$ 
\begin{align*}
 u(s,X_s^{t,x})&=
 \E_s[\phi(X_T^{t,x})]+\E_s\left[\int_s^T
 \psi(\tau,X_\tau^{t,x}, u(\tau,X_\tau^{t,x} ),
  \nabla u(\tau,X_\tau^{t,x})B)\,d\tau\right]\\ \nonumber
 &=\E_s[\xi]-\int_t^s
 \psi(\tau,X_\tau^{t,x}, u(\tau,X_\tau^{t,x} ), \nabla u(\tau,X_\tau^{t,x})B)\,d\tau,\\ \nonumber
\end{align*}
where
\[
\xi:=\phi(X_T^{t,x})+\int_t^T
\psi(\tau,X_\tau^{t,x}, u(\tau,X_\tau^{t,x} ), \nabla u(\tau,X_\tau^{t,x})B)\,d\tau.
\]
By the martingale representation theorem,
there exists a process $\tilde Z\in\,L^2(\Omega\times[t,T];\Xi^{\ast})$
such that $\E_s[\xi]=u(t,x)+\int_t^s\tilde Z_\tau\,dW_\tau$.
So $ \left(u(s,X_s^{t,x})\right)_{s\in[t,T]}$ is a continuous semi-martingale
with canonical decomposition
\begin{equation}\label{u-semimart}
 u(s,X_s^{t,x})=u(t,x)+\int_t^s\tilde Z_\tau\,dW_\tau-\int_t^s
\psi(\tau,X_\tau^{t,x}, u(\tau,X_\tau^{t,x} ), \nabla u(\tau,X_\tau^{t,x})B)\,d\tau.
\end{equation}
As in the Lemma 6.3 of \cite{fute}, when we compute the joint quadratic variation
of $ \left(u(s,X_s^{t,x})\right)_{s\in[t,T]}$ with the Wiener process, we get the identification
\[
  \nabla u(\tau,X_\tau^{t,x})B=\tilde Z _\tau^{t,x}.
\]
Substituting into (\ref{u-semimart}), and rewriting the obtained equality
in backward sense, we note that $\left(Y^{t,x},Z^{t,x}\right)$ and $\left(u(\cdot, X^{t,x}) ,\nabla u(\cdot, X^{t,x})B\right)$ solve the same equation, and so uniqueness follows from the uniqueness of the BSDE solution.
\qed

\section{Mild solution of a semilinear PDE in infinite dimension: the Lipschitz
continuous data case}
\label{sezionePDElip}
The aim of this section is to study equation (\ref{Kolmo})
when the final datum $\phi$ and the nonlinear term $\psi$
are only Lipschitz continuous.
Notice that in order to do this, we require some smoothing
properties on the transition semigroup $(P_{t,\tau})_{\tau\in[t,T]}$. Namely we require the following smoothing property on the semigroup
$(P_{t,\tau})_{\tau\in[t,T]}$, see e.g. \cite{Mas} where this property has been introduced for bounded functions,
and \cite{Mas2} where it has been extended to functions with polynomial growth.
\begin{hypothesis}\label{ipH su fi} For some $\alpha\in [0,1)$ and for
every $\phi\in C_{k}\left(  H\right)  $, the function $P_{t,\tau}\left[
\phi\right]  \left(  x\right)  $ is $B  $-differentiable
with respect to $x$, for every $0\leq t <\tau < T$. Moreover, for every $k \in \mathbb{N}$ there exists a
constant $c_k>0$ such that for every $\phi\in C_{k}\left(  H\right)  $, for
every $\xi\in\Xi,$ and for $0\leq t<\tau\leq T$,%
\begin{equation*}
\left|  \nabla^BP_{t,\tau}\left[  \phi\right]
\left(  x\right)  \xi\right|  \leq\frac{c_k}{\left(\tau-t\right)  ^{\alpha}%
}\left\|  \phi\right\|  _{C_k}\left|  \xi\right|  .
\end{equation*}
\end{hypothesis}
In \cite{Mas} it is shown that Hypothesis \ref{ipH su fi} is verified for Ornstein-Uhlenbeck transition
semigroups ( i.e. $F=0$ in (\ref{ornstein-pert})) by relating $B$-differentiability
to properties of the operators $A$ and $B$, as collected in the following proposition.
\begin{proposition}
 Let us assume that 
 \begin{equation}
\operatorname{Im}e^{(\tau-t)A}B\subset\operatorname{Im}Q_{\tau-t}^{1/2},
\label{ornstein inclusione}%
\end{equation}
and, for some $0\leq\alpha<1$ and $c>0$, the operator norm satisfies
\begin{equation}
 \left\|  Q_{\tau-t}^{-1/2}e^{(\tau-t)A}B\right\|  \leq c(\tau-t)^{-\alpha} \quad \text{ for } \quad
0\leq t<\tau\leq T.\label{ornstein inclusionestima}
\end{equation}
Then Hypothesis \ref{ipH su fi} is satisfied by the Ornstein-Uhlenbeck transition semigroup.
\end{proposition}

We refer to \cite{Mas} where some examples of Ornstein-Uhlenbeck semigroup satisfying Hypothesis
\ref{ipH su fi} are provided. Among these examples we remember the wave equation, see
also section \ref{sez-contr-wave}.

\noindent We now prove existence and uniqueness of a mild solution for the Kolmogorov equation (\ref{Kolmo})
when $\mathcal L$ is the generator of a Ornstein-Uhlenbeck transition semigroup, that is to say $F=0$ in (\ref{ornstein-pert}). 
The perturbed Ornstein-Uhlenbeck case will be treated after in Theorem \ref{teokolmolippert}.

\begin{theorem}\label{teokolmolip}
Assume that Hypotheses \ref{ip su AB} and \ref{ip-psi-phi}
hold true, and let $F=0$ in (\ref{ornstein-pert}), and consequently also in (\ref{Kolmo}) so that the process $X^{t,x}$ is an Ornstein-Uhlenbeck process.
Moreover assume that the Ornstein-Uhlenbeck transition semigroup related to $X^{t,x}$
 satisfies Hypothesis \ref{ipH su fi}. Then, according to definition \ref{defsolmildkolmo},
equation (\ref{Kolmo}) admits a unique mild solution.
\end{theorem}

\textbf{Proof.} The idea of the proof is to 
smooth coefficients $\psi $ and $\phi$, so
to obtain a sequence of approximating Kolmogorov equations which admit a solution
according to Theorem \ref{teokolmosmooth}, and then to pass to the limit.

\noindent Coming into more details, we are approximating functions that 
have polynomial growth with respect to their arguments and are (locally) Lipschitz continuous, but we
need to preserve their (locally) Lipschitz constant. So to approximate these functions we follow
\cite{PZ}. In that paper for every
$n\in\mathbb{N}$ it is considered a nonnegative function $\rho_{n}\in C_{b}^{\infty}\left(  \mathbb{R}%
^{n}\right)  $ with compact support contained in the
ball of radius $\frac{1}{n}$ and such that $
{\displaystyle\int_{\mathbb{R}^{n}}}
\rho_{n}\left(  x\right)  dx=1$. Let $\left\{  e_{k}\right\}  _{k\in
\mathbb{N}}$ be a complete orthonormal system in $H$ and, for every
$n\in\mathbb{N}$, let $Q_{n}:H\longrightarrow\left\langle e_{1},...,e_{n}%
\right\rangle $ be the orthogonal projection on the linear space generated by
$e_{1},...,e_{n}$. We identify $\left\langle e_{1},...,e_{n}\right\rangle $
with $\mathbb{R}^{n}$. For a bounded and continuous function $f:H\rightarrow \R$ we set
\[
f_{n}\left(  x\right)  =\int_{\mathbb{R}^{n}}\rho_{n}\left(  y-Q_{n}%
x\right)  f\left(  \sum_{i=1}^{n}y_{i}e_{i}\right)  dy,
\]
where for every $k\in\mathbb{N}$, $y_{k}=\left\langle y,e_{k}\right\rangle
_{H}$. It turns out that $f_{n}\in C_{b}^{\infty}\left(  H\right)  $.
Moreover, if $f$ is (locally) Lipschitz continuous and has polynomial growth, $f_n$ is
(locally) Lipschitz continuous and has polynomial growth as well, it
preserves the (locally) Lipschitz constant and the order of polynomial growth is the same as the one of $f$.
Namely, if there exist $L>0$ and $C>0$ such that
\[
\left\vert f\left(  x\right)  -f\left(  y\right)  \right\vert \leq
L\left\vert x-y\right\vert\left(  1+\left\vert
x\right\vert^r+\vert y\vert^r \right) ,\text{ \ \ for every }x,y\in H,
\]
then for every $k\in\mathbb{N}$%
\[%
\left\vert f_{n}\left(  x\right)  -f_{n}\left(  y\right)  \right\vert
\leq L\left\vert x-y\right\vert\left(  1+\left\vert
x\right\vert^r+\vert y\vert^r \right) ,\text{ \ \ for every }x,y\in H.
\]
Finally, $\left(  f_{n}\right)  _{n}$ is a pointwise approximation of $f$: for every $x\in H,$%
\[
\lim_{n\rightarrow\infty}\left\vert f_{n}\left(  x\right)  -f\left(
x\right)  \right\vert =0.%
\]
So, if we consider the final datum $\phi$ in the Kolmogorov equation \ref{Kolmo},
we can set
\begin{equation}\label{infsupconvphi}
\phi_{n}\left(  x\right)  =\int_{\mathbb{R}^{n}}\rho_{n}\left(  y-Q_{n}%
x\right)  \phi\left(  \sum_{i=1}^{n}y_{i}e_{i}\right)  dy,
\end{equation}
and we have that, $\forall x,\, x'\,\in H$ and $n\in\N$
\[
 \vert \phi_n(x)-\phi_n(x')\vert \leq 
\left(C+\dfrac{\alpha}{2}\vert x\vert^r+\dfrac{\alpha}{2}\vert x'\vert^r\right)\vert x-x'\vert.
\]
For what concerns $\psi$, we consider another sequence of functions $(\bar \rho _n)_n$
satisfying the same properties introduced before for the sequence $(\rho_n)_n$, and $\left\{   \bar e _{k}\right\}  _{k\in \N} $ 
a complete orthonormal system in $\Xi^*$. Finally let $(\hat{\rho}_n)_n$ be a sequence of nonnegative real functions with compact support 
contained in $[-1/n,1/n]$ and such that ${\displaystyle\int_{\mathbb{R}}} \hat{\rho}_{n}\left(  x\right)  dx=1$.
So we can define
\begin{equation}\label{infsupconvpsi}
 \psi_{n}\left(t,  x,y,z\right)  =\int_{\mathbb{R}^{n}}\int_\R\int_{\R^n}\rho_{n}\left(  x'-Q_{n}%
x\right) \hat\rho_n (y'-y) \bar\rho_n (z'-\bar{Q}_n z)\psi\left( t, \sum_{i=1}^{n}x'_{i}e_{i},y',\sum_{i=1}^{n}z'_{i}\bar e_{i}\right)
 dx'\, dy'\,dz'.
\end{equation}
We have that for all $t\in [0,T],\,x,\,x'\in H,\,y,y'\in \R,\,z,z'\in \Xi^*$,
\begin{align*}
&\psi_n(t,x,y,z)-\psi(t,x,y',z)\vert \leq K_{\psi_y}\vert y -y'\vert;\\
& \vert \psi_n(t,x,y,z)-\psi_n(t,x,y,z')\vert \leq 
\left(C+\frac{\gamma}{2}\vert z\vert^l+\frac{\gamma}{2}\vert z'\vert^l\right)\vert z -z'\vert;\\
&  \vert \psi_n(t,x,y,z)-\psi_n(t,x',y,z)\vert \leq 
\left(C+\frac{\beta}{2}\vert x\vert^r+\frac{\beta}{2}\vert x'\vert^r\right)\vert x -x'\vert.
\end{align*}

We notice that we only have a pointwise convergence of $\phi_n$ to $\phi$ and of $\psi_n$ to $\psi$, see again \cite{PZ}.
For this reason in the sequel it will be crucial the fact that $P$ is an Ornstein-Uhlenbeck transition semigroup, so that we can explicitly represent the mild solution of the Kolmogorov equation.

Now the proof goes on by approximating $\phi$ and $\psi$,
so to build a sequence of mild solutions of the Kolmogorov equations with the approximating coefficients
$\phi_n$ and $\psi_n$. We want to prove that the
sequence of solutions converges in a suitable space. 
Firstly, we need a stability result for the solution of the BSDE
(\ref{bsde}) with respect to the approximation of the final datum and the generator.
\begin{proposition}
 \label{prop-approx-bsde}
Let $(Y^{n,t,x},Z^{n,t,x})$ and $(Y^{k,t,x},Z^{k,t,x})$ be solutions of the BSDE (\ref{bsde})
with final datum and generator respectively given by the approximants
$\phi_n$ and $\psi_n$, and by $\phi_k$ and $\psi_k$ defined respectively in (\ref{infsupconvphi}) and in (\ref{infsupconvpsi}).
Namely 
\begin{align}
Y^{n,t,x}_\tau-Y^{k,t,x}_\tau=&\phi_n(X_T^{t,x})-\phi_k(X_T^{t,x})-\int_\tau^T\left(Z^{n,t,x}_s-Z^{k,t,x}_s\right)\,dW_s \nonumber \\
&+\int_\tau^T \left(\psi_n(s,X_s^{t,x}, Y^{n,t,x}_s,Z^{n,t,x}_s)
-\psi_k(s,X_s^{t,x}, Y^{k,t,x}_s,Z^{n,t,x}_s)\right)\,ds. \nonumber
\end{align}
Then, $\forall\, t\in[0,T],\,x \in H$, we have
\begin{equation*}
 \Vert Y^{n,t,x}-Y^{k,t,x}\Vert_{\cals^2}+\Vert Z^{n,t,x}-Z^{k,t,x}\Vert_{\calm^2} \leq C_{n,k}(t,x),
\end{equation*}
with $\displaystyle{\lim_{n,k\rightarrow \infty} C_{n,k}(t,x)= 0}$.
\end{proposition}
 \textbf{Proof of Proposition \ref{prop-approx-bsde}.} 
By the usual linearization trick we can write 
\begin{align*}
&Y^{n,t,x}_\tau-Y^{k,t,x}_\tau=\phi_n(X_T^{t,x})-\phi_k(X_T^{t,x})+\int_\tau^T \left(\psi_n(s,X_s^{t,x}, Y^{n,t,x}_s,Z^{n,t,x}_s)
-\psi_k(s,X_s^{t,x}, Y^{n,t,x}_s,Z^{n,t,x}_s)\right)\,ds\\
&+\int_\tau^T U^{n,k}_s\left( Y^{n,t,x}_s-Y^{k,t,x}_s\right)\,ds
+\int_\tau^T \langle V^{n,k}_s,\left( Z^{n,t,x}_s-Z^{k,t,x}_s\right)\rangle\,ds-\int_\tau^T\left(Z^{n,t,x}_s-Z^{k,t,x}_s\right)\,dW_s\
\end{align*}
where we have set
\[
U^{n,k}_s=\left\lbrace
\begin{array}{ll}
\dfrac{\psi_k(s,X_s^{t,x}, Y^{n,t,x}_s,Z^{n,t,x}_s)
-\psi_k(s,X_s^{t,x}, Y^{k,t,x}_s,Z^{n,t,x}_s) }{Y^{n,t,x}_s-Y^{k,t,x}_s} &\text{if }Y^{n,t,x}_s-Y^{k,t,x}_s\neq 0\\
0&\text{if }Y^{n,t,x}_s-Y^{k,t,x}_s= 0,
\end{array}
 \right.\]
and\[
 V^{n,k}_s=\left\lbrace
 \begin{array}{ll}
 \dfrac{\psi_k(s,X_s^{t,x}, Y^{k,t,x}_s,Z^{n,t,x}_s)
 -\psi_k(s,X_s^{t,x}, Y^{k,t,x}_s,Z^{k,t,x}_s) }{\abs{Z^{n,t,x}_s-Z^{k,t,x}_s}^2}(Z^{n,t,x}_s-Z^{k,t,x}_s) &\text{if }Z^{n,t,x}_s-Z^{k,t,x}_s\neq 0\\
 0&\text{if }Z^{n,t,x}_s-Z^{k,t,x}_s= 0.
 \end{array}
 \right.
   \]
Since $\vert V^{n,k}_s\vert\leq C\left(1+\vert X_s^{t,x}\vert^{rl}\right)$, by the Girsanov theorem there exists a probability measure $\Q^{n,k}$, equivalent to the original one $\P$, such that
$\tilde W _\tau:=W_\tau-\int_0^\tau V^{n,k}_s\,ds$ is a $\Q^{n,k}$-Wiener process and we have
\begin{align*}
 Y^{n,t,x}_\tau-Y^{k,t,x}_\tau&=\E^{\Q^{n,k}}_\tau\left[e^{\int_\tau^TU^{n,k}_s\,ds}
\left(\phi_n(X_T^{t,x})-\phi_k(X_T^{t,x})\right)\right]\\
&+\E^{\Q^{n,k}}_\tau\left[\int_t^T e^{\int_s^TU^{n,k}_r\,dr}\left(\psi_n(s,X_s^{t,x}, Y^{n,t,x}_s,Z^{n,t,x}_s)
-\psi_k(s,X_s^{t,x}, Y^{n,t,x}_s,Z^{n,t,x}_s)\right)\,ds\right].
\end{align*}
Since $\vert U^{n,k}_s\vert\leq CK_\psi$, we get
 \begin{align*}
 \vert Y^{n,t,x}_\tau-Y^{k,t,x}_\tau\vert^2&\leq C\E^{\Q^{n,k}}_\tau \left[\abs{\phi_n\left (X^{t,x}_T\right )-
 \phi_k\left (X^{t,x}_T\right)}\right]^2\\
 &+C\E^{\Q^{n,k}}_\tau \left[\int_{t}^T \vert \psi_n(s,X_s^{t,x}, Y^{n,t,x}_s,Z^{n,t,x}_s)
 -\psi_k(s,X_s^{t,x}, Y^{n,t,x}_s,Z^{n,t,x}_s)\vert ds\right]^2.
 \end{align*}
By keeping in mind that $\vert Z_s^{n,t,x}\vert \leq C\left(1+\vert X_s^{t,x}\vert^r\right)$
and $\vert Y_s^{n,t,x}\vert \leq C\left(1+\vert X_s^{t,x}\vert^{r+1}\right)$, we have
$$\vert \psi_n(s,X_s^{t,x}, Y^{n,t,x}_s,Z^{n,t,x}_s)-\psi_k(s,X_s^{t,x}, Y^{n,t,x}_s,Z^{n,t,x}_s) \vert \leq C\left(1+\vert X^{t,x}_s\vert^{r+1} \right), $$
and the dominated convergence theorem gives us
\begin{equation*}
\E\left[\sup_{\tau\in[0,T]} \vert Y^{n,t,x}_\tau-Y^{k,t,x}_\tau\vert^2\right]\rightarrow
0 \quad\text{as }n,k\rightarrow \infty,
\end{equation*}
pointwise with respect to $t$ and $x$.
Now we look for an estimate for the $\calm^2$-norm of $Z^{n,t,x}-Z^{k,t,x}$.
By applying It\^o formula to $\vert Y^{n,t,x}_\tau-Y^{k,t,x}_\tau\vert^2$
we get
\begin{align*}
\vert Y^{n,t,x}_0-Y^{k,t,x}_0\vert^2=&\vert\phi_n(X_T^{t,x})-\phi_k(X_T^{t,x})\vert^2\\
&+2\int_0^T \left(Y^{n,t,x}_s-Y^{k,t,x}_s\right) \left(\psi_n(s,X_s^{t,x}, Y^{n,t,x}_s,Z^{n,t,x}_s)
-\psi_k(s,X_s^{t,x}, Y^{k,t,x}_s,Z^{k,t,x}_s)\right)\,ds\\
&-\int_0^T\abs{ Z^{n,t,x}_s-Z^{k,t,x}_s}^2\,ds-2\int_0^T\left( Y^{n,t,x}_s-Y^{k,t,x}_s\right)\left(Z^{n,t,x}_s-Z^{k,t,x}_s\right)\,dW_s.
\end{align*}
By taking expectation and by standard calculations we get
\begin{align*}
\E \int_0^T\abs{ Z^{n,t,x}_s-Z^{k,t,x}_s}^2\,ds
\leq&\E \vert\phi_n(X_T^{t,x})-\phi_k(X_T^{t,x})\vert^2\\
&+C\E\int_\tau^T \abs{Y^{n,t,x}_s-Y^{n,t,x}_s} \left(1+\abs{X_s^{t,x}}^{r+1}\right)\,ds.
 \end{align*}
So, $\forall\, t\in[0,T], \, x\in H$, $\Vert Z^{n,t,x}-Z^{k,t,x}\Vert_{\calm^2}\rightarrow 0$ as $n,k \rightarrow \infty$, and the proposition is proved.
\qed

\noindent Next we go on proving Theorem \ref{teokolmolip}.

\noindent\textbf{Proof of Theorem \ref{teokolmolip}-continuation.}
We denote by $v^n$ the solution of the Kolmogorov equation (\ref{Kolmo}), with final datum
$\phi_n$ instead of $\phi$ and Hamiltonian function $\psi_n$ instead of $\psi$. Namely
$v_n$ satisfies
\begin{equation}
v_n(t,x)=P_{t,T}\left[  \phi_n\right]  \left(  x\right)  +\int_{t}^{T}%
P_{t,s}\left[  \psi_n\left(s,\cdot,v_n\left(  s,\cdot\right),\nabla^B 
v_n\left(  s,\cdot\right)\right)  \right] (  x)  ds.
 \label{solmildkolmon}%
\end{equation}
Since the data $\phi_n$ and $\psi_n$ are differentiable, we also know by theorem \ref{teokolmosmooth}
that the pair of processes $(v_n(\cdot, X^{t,x}_{\cdot}),\nabla^Bv_n(\cdot, X^{t,x}_{\cdot}))$ is solution
to the following BSDE
\begin{equation*}
  Y_s^{n,t,x}+\int_s^T Z^{n,t,x}_\tau\;dW_\tau=\phi_n(X_T^{t,x})+\int_s^T\psi_n(\tau, X^{t,x}_\tau,Y_\tau^{n,t,x},Z_\tau^{n,t,x})\;d\tau,
\end{equation*}
so we get that, for every $n\in\N$, $t\in[0,T]$ and $x\in H$
\begin{equation*}
 \vert v_n(t,x) \vert \leq C\left(1+\vert x\vert^{r+1} \right), \qquad
\vert \nabla^B v_n(t,x) \vert \leq C\left(1+\vert x\vert^{r} \right),
\end{equation*}
where $C$ is a constant that does not depend on $n$, $t$, $x$, see Proposition \ref{prop-exist-bsde}.
We want to show that $v_n$ converges to $v$, a solution of the Kolmogorov equation (\ref{Kolmo}).
By Proposition \ref{prop-approx-bsde} we know that the sequence $(v_n(t,x))_{n\geq1}$ is a Cauchy sequence for all $t \in [0,T]$, $x \in H$,
and we want to show that the sequence $(\nabla^B v_n(t,x))_{n\geq1}$ is a Cauchy sequence for all $t \in [0,T[$, $x \in H$.
Let us recall that, by identification (\ref{identif-Z}) of $Z$, we have
\[ \vert\nabla^B v_n(t,x)\vert \leq C(1+\vert x\vert ^2)^{r/2},
\]
with $C$ a constant independent on $n$, $t$, $x$. Notice that, in virtue of Hypothesis \ref{ip-psi-phi}, and of this estimate,
the map $x\mapsto \psi (s,x,v_n(s,x),\nabla^B v_n(s,x))$ has polynomial growth
of order $r+1$ uniformly with respect to $s\in[t,T]$ and to $n\geq 1$, that is 
\[
 \vert \psi (s,x,v_n(s,x),\nabla^B v_n(s,x))\vert \leq C \left( 1+\vert x\vert ^2 \right)^{(r+1)/2},
\]
with $C$ a constant independent on $n,\, s$ and $x$.

We consider, for $n,k\geq 1$, the difference $ v_n(t,x)-v_k(t,x)$
\begin{align*}
 v_n(t,x)-v_k(t,x)=& P_{t,T}\left[  \phi_n-\phi_k\right]  \left(  x\right) \\
&+\int_{t}^{T}P_{t,s}\left[
\psi_n(s,\cdot, v_n(s,\cdot),\nabla^Bv_n(s,\cdot))-
\psi_k(s,\cdot, v_k(s,\cdot),\nabla^Bv_k(s,\cdot))\right](x)\,ds\\
=&\int_{H}
\left(\phi_n(z+e^{(T-t)A}x)-\phi_k(z+e^{(T-t)A}x)\right)
\caln(0,Q_{s-t})(dz)\\
&+\int_{t}^{T}\int_{H}\left[
\psi_n\left(s,x,v_n (s,z+e^{(s-t)A}x),\nabla^Bv_n(s,z+e^{(s-t)A}x)\right)\right.\\
&\left.-\psi_k\left(s,x,v_k (s,z+e^{(s-t)A}x),\nabla^Bv_k(s,z+e^{(s-t)A}x)\right)\right]
\caln(0,Q_{s-t})(dz).
\end{align*}
Since $v_n$ and $v_k$ are G\^ateaux differentiable and by the smoothing properties of
the transition semigroup $(P_{t,\tau})_{\tau\in[t,T]}$, we can take the $B$ derivative of both sides in
(\ref{solmildkolmon}) and, by the closedness of the operator $\nabla^B$, see e.g. \cite{Mas},
we obtain for all $h \in \Xi$
\begin{align*}
&\nabla^B v_n(t,x)h-\nabla^B v_k(t,x)h =\nabla ^B P_{t,T}\left[  \phi_n-\phi_k\right]  \left(  x\right)h \\
 &+\int_{t}^{T}\nabla^BP_{t,s}\left[  \psi_n\left(s,\cdot,v_n\left(  s,\cdot\right),\nabla^B 
v_n\left(  s,\cdot\right)\right)-\psi_k\left(s,\cdot,v_k\left(  s,\cdot\right),\nabla^B 
v_k\left(  s,\cdot\right)\right)  \right] (  x)h  ds.
\end{align*}
Namely, following \cite{Mas},
when $X$ is an Ornstein-Uhlenbeck process we have an explicit expression for the $B$-derivative of the transition semigroup applied to some continuous function,
see Lemma 3.4 in \cite{Mas}, generalized to the case of functions with polynomial growth with respect to $x$ in \cite{Mas2}.
We get that for every continuous function $f\in C_{r+1}(H)$ and every $h \in \Xi$ we have
\begin{equation*}
 \nabla^{B}(P_{t,s}\left[  f\right]  )\left(  x\right)  h\\
 =\int_{H}f\left(  y+e^{(s-t)A}x\right)  \left\langle Q_{s-t}^{-1/2}%
e^{(s-t)A}Bh,Q_{s-t}^{-1/2}y\right\rangle \mathcal{N}\left(  0,Q_{s-t}\right)  \left(
dy\right)  .
\end{equation*}
Taking into account this fact, we get
\begin{eqnarray*}
&&\nabla^B v_n(t,x)h-\nabla^B v_k(t,x)h \\
&=&\int_{H}\left[\phi_n\left(  z+e^{(T-t)A}x\right)-\phi_k\left(  z+e^{(T-t)A}x\right)\right]
  \left\langle Q_{T-t}^{-1/2}%
e^{(T-t)A}Bh,Q_{T-t}^{-1/2}y\right\rangle \mathcal{N}\left(  0,Q_{T-t}\right)  \left(
dz\right)\\
&&+\int_t^T\int_{H}\left[\psi_n\left(s,  z+e^{(s-t)A}x, v_n(s,  z+e^{(s-t)A}x), \nabla^Bv_n(s,  z+e^{(s-t)A}x)\right)\right.\\
&&\left.\;-\psi_k\left(s,  z+e^{(s-t)A}x, v_k(s,  z+e^{(s-t)A}x),\nabla^Bv_k(s,  z+e^{(s-t)A}x)\right)\right]\\
 & &\times \left\langle Q_{s-t}^{-1/2}%
e^{(s-t)A}Bh,Q_{s-t}^{-1/2}y\right\rangle \mathcal{N}\left(  0,Q_{s-t}\right)  \left(
dz\right)ds.
\end{eqnarray*}
Now we want to estimate $\vert \nabla^B v_n(t,x)h-\nabla^B v_k(t,x)h \vert$.
At first we consider
\begin{eqnarray*}
& &\vert \nabla^B P_{t,T}\left[  \phi_n-\phi_k\right]  \left(  x\right)h\vert\\
&=& \abs{\int_H \left(\phi_n\left(  z+e^{(T-t)A}x\right)-\phi_k\left(  z+e^{(T-t)A}x\right)\right)
 \left\langle Q_{T-t}^{-1/2}%
e^{(T-t)A}Bh,Q_{T-t}^{-1/2}y\right\rangle \mathcal{N}\left(  0,Q_{T-t}\right)  \left(
dz\right)}\\
&  \leq &\left( \int_H \abs{ \phi_n\left(  z+e^{(T-t)A}x\right)-\phi_k\left(  z+e^{(T-t)A}x\right)}^2\mathcal{N}\left(  0,Q_{T-t}\right)  \left(
 dz\right)
  \right)^{1/2}\\
&&\times \left( \int_H \abs{  \left\langle Q_{T-t}^{-1/2}%
 e^{(T-t)A}Bh,Q_{T-t}^{-1/2}y\right\rangle } ^2\mathcal{N}\left(  0,Q_{T-t}\right)  \left(
 dz\right)  \right)^{1/2} \\
  &\leq& C \left(T-t\right)^{-1/2} \left(\int_H \abs{ \phi_n\left(  z+e^{(T-t)A}x\right)-\phi_k\left(  z+e^{(T-t)A}x\right)}^2\mathcal{N}\left(  0,Q_{T-t}\right)  \left(
  dz\right)
   \right)^{1/2}\abs{h},
\end{eqnarray*}
and so $\vert \nabla^B P_{t,T}\left[  \phi_n-\phi_k\right]  \left(  x\right)h\vert$
converges pointwise to $0$ for all $x\in H$ and $t\in [0,T)$ as $n,k \rightarrow + \infty$.
Now we have to estimate 
\begin{align*}
&\int_{t}^{T}%
 \nabla^BP_{t,s}\left[  \psi_n\left(s,\cdot,v_n\left(  s,\cdot\right),\nabla^B 
 v_n\left(  s,\cdot\right)\right)-\psi_k\left(s,\cdot,v_k\left(  s,\cdot\right),\nabla^B 
 v_k\left(  s,\cdot\right)  \right)\right] (  x)  ds\\
  =&\int_{t}^{T}\nabla^BP_{t,s}\left[  \psi_n\left(s,\cdot,v_n\left(  s,\cdot\right),\nabla^B 
  v_n\left(  s,\cdot\right)\right)-\psi_k\left(s,\cdot,v_n\left(  s,\cdot\right),\nabla^B 
  v_n\left(  s,\cdot\right) \right) \right] (  x) h \, ds\\
  &+\int  \nabla^BP_{t,s}\left[  \psi_k\left(s,\cdot,v_n\left(  s,\cdot\right),\nabla^B 
  v_n\left(  s,\cdot\right)\right)-\psi_k\left(s,\cdot,v_k\left(  s,\cdot\right),\nabla^B 
  v_n\left(  s,\cdot\right)\right)  \right] (  x)  h\, ds\\
  &+\int_{t}^{T}\nabla^BP_{t,s}\left[  \psi_k\left(s,\cdot,v_k\left(  s,\cdot\right),\nabla^B 
  v_n\left(  s,\cdot\right)\right)-\psi_k\left(s,\cdot,v_k\left(  s,\cdot\right),\nabla^B 
  v_k\left(  s,\cdot\right)\right)  \right] (  x) h\, ds\\
=&I+II+III.
\end{align*}
With calculations similar to the ones performed for estimating 
\[
\vert P_{t,T}\left[  \phi_n-\phi_k\right]  \left(  x\right)h\vert +\vert \nabla^B P_{t,T}\left[  \phi_n-\phi_k\right] 
\left(  x\right)h\vert
\]
we get
\begin{align*}
\vert I\vert=&\bigg\vert\int_t^T\int_H \left[\psi_n\left( s, z+e^{(s-t)A}x, v_n(s, z+e^{(s-t)A}x),\nabla^B v_n(s, z+e^{(s-t)A}x)\right)\right.\\
&\left.\;-\psi_k\left( s, z+e^{(s-t)A}x, v_n(s, z+e^{(s-t)A}x),\nabla^B v_n(s, z+e^{(s-t)A}x)\right)\right]\\
 &\times\left\langle Q_{T-t}^{-1/2}%
 e^{(T-t)A}Bh,Q_{T-t}^{-1/2}y
\right\rangle \mathcal{N}\left(  0,Q_{T-t}\right)  \left(
 dz\right)\,ds\bigg\vert\\
 \leq &C \left( T-t\right)^{-1/2} \left(\int_t^T\int_H \bigg\vert\psi_n\left( s, z+e^{(s-t)A}x, v_n(s, z+e^{(s-t)A}x),\nabla^B v_n(s, z+e^{(s-t)A}x)\right)\right.\\
 &\left.-\psi_k\left( s, z+e^{(s-t)A}x, v_n(s, z+e^{(s-t)A}x),\nabla^B v_n(s, z+e^{(s-t)A}x)\right)\bigg\vert^2\mathcal{N}\left(  0,Q_{T-t}\right)  \left(
   dz\right)ds
    \right)^{1/2}\abs{h}\\
    & \rightarrow 0 \qquad \text{as}\;n,k\rightarrow \infty,
\end{align*}
pointwise for all $x\in H$ and $t\in [0,T)$, by the dominated convergence theorem and by the convergence
of $\psi_n$, as well of $\psi_k$, to $\psi$.
Next we estimate $II$:
\begin{align*}
 \vert II\vert=&\bigg\vert\int_t^T\int_{H}\left[\psi_k\left(s,  y+e^{(s-t)A}x,v_n(s,y+e^{(s-t)A}x),\nabla^Bv_n( y+e^{(s-t)A}x)\right)  \right.\\
&\left.-\psi_k\left(s,  y+e^{(s-t)A}x,v_k(s,y+e^{(s-t)A}x),\nabla^Bv_n( y+e^{(s-t)A}x)\right)\right]\\
&\times \left\langle Q_{s-t}^{-1/2} e^{(s-t)A}Bh,Q_{s-t}^{-1/2}y\right\rangle \mathcal{N}\left(  0,Q_{s-t}\right)(dy)\,ds\bigg\vert\\
 \leq& \int_t^T\left(\int_{H}\Big\vert\psi_k\left(s,  y+e^{(s-t)A}x,v_n(s,y+e^{(s-t)A}x),\nabla^Bv_n( y+e^{(s-t)A}x)\right)  \right.\\
&-\psi_k\left(s,  y+e^{(s-t)A}x,v_k(s,y+e^{(s-t)A}x),\nabla^Bv_n( y+e^{(s-t)A}x)\right)\Big\vert^2\mathcal{N}\left(  0,Q_{s-t}\right)(dy)\bigg)^{1/2}\\
& \times \left(\int_H\bigg\vert\left\langle Q_{s-t}^{-1/2} e^{(s-t)A}Bh,Q_{s-t}^{-1/2}y\right\rangle\bigg\vert^2 \mathcal{N}\left(  0,Q_{s-t}\right)(dy)\right)^{1/2}\,ds\\
 \leq& C\int_t^T(s-t)^{-\alpha}\left(\int_{H}\vert v_n( y+e^{(s-t)A}x)  - v_k( y+e^{(s-t)A}x)\vert^2\mathcal{N}\left(  0,Q_{s-t}\right)(dy)\right)^{1/2}ds\abs{h}\\
 & \rightarrow 0 \qquad \text{as}\;n,k\rightarrow \infty
 \end{align*}
for all $t, x\in\, H$, where in the last passage we have used the dominated convergence theorem and the
pointwise convergence of $v_n-v_k$ to $0$.
Finally we estimate $III$:
\begin{align*}
 \vert III\vert=&\bigg\vert\int_t^T\int_{H}\left[\psi_k\left(s,  y+e^{(s-t)A}x,v_k(s,y+e^{(s-t)A}x),\nabla^Bv_n( y+e^{(s-t)A}x)\right)  \right.\\
&\left.-\psi_k\left(s,  y+e^{(s-t)A}x,v_k(s,y+e^{(s-t)A}x),\nabla^Bv_k( y+e^{(s-t)A}x)\right)\right]\\
&\times \left\langle Q_{s-t}^{-1/2} e^{(s-t)A}Bh,Q_{s-t}^{-1/2}y\right\rangle \mathcal{N}\left(  0,Q_{s-t}\right)(dy)\,ds\bigg\vert\\
 \leq& \int_t^T\left(\int_{H}\Big\vert\psi_k\left(s,  y+e^{(s-t)A}x,v_n(s,y+e^{(s-t)A}x),\nabla^Bv_n( y+e^{(s-t)A}x)\right)  \right.\\
&-\psi_k\left(s,  y+e^{(s-t)A}x,v_k(s,y+e^{(s-t)A}x),\nabla^Bv_n( y+e^{(s-t)A}x)\right)\Big\vert^2\mathcal{N}\left(  0,Q_{s-t}\right)(dy)\bigg)^{1/2}\\
& \times \left(\int_H\bigg\vert\left\langle Q_{s-t}^{-1/2} e^{(s-t)A}Bh,Q_{s-t}^{-1/2}y\right\rangle\bigg\vert^2 \mathcal{N}\left(  0,Q_{s-t}\right)(dy)\right)^{1/2}\,ds\\
 \leq& C\int_t^T(s-t)^{-\alpha}\left(\int_{H}\vert\nabla^Bv_n( y+e^{(s-t)A}x)  - \nabla^Bv_k( y+e^{(s-t)A}x)\vert^2\right.\\
&\times\left(1+\vert\nabla^Bv_n(s,y+e^{(s-t)A}x)\vert^{2l}+\vert\nabla^Bv_k(s,y+e^{(s-t)A}x)\vert^{2l}\right)\mathcal{N}\left(  0,Q_{s-t}\right)(dy)\bigg)^{1/2}ds\abs{h}\\
 =&C\int_t^T(s-t)^{-\alpha}\left( \E\left[\vert Z^{n,t,x}_s-Z^{k,t,x}_s\vert^2
 \left(1+\vert Z^{n,t,x}_s\vert^{2l}+\vert Z^{k,t,x}_s\vert^{2l}\right)\right]\right)^{1/2}\,ds\abs{h}.
 \end{align*}

 Then, by using the uniform bound (with respect to $n$) on $Z^{n,t,x}$ and $Z^{k,t,x}$ and by the H\"older inequality, we obtain
 \begin{align*}
 \vert III\vert\leq& C\int_t^T(s-t)^{-\alpha}\left( \E\left[\vert Z^{n,t,x}_s-Z^{k,t,x}_s\vert^{1-\alpha}
 \left(1+\vert X^{t,x}_s\vert^{2rl+r+r\alpha}\right)\right]\right)^{1/2}\,ds\abs{h}\\
 \leq& C\int_t^T(s-t)^{-\alpha} \E\left[\vert Z^{n,t,x}_s-Z^{k,t,x}_s\vert^2\right]^{(1-\alpha)/4}
 \E\left[1+\vert X^{t,x}_s\vert^{2(2rl+r+r\alpha)/(1+\alpha)}\right]^{(1+\alpha)/4}\,ds\abs{h}\\
\leq& C\sup_{s\in[0,T]} \left(\E\left[1+\vert X^{t,x}_s\vert^{2(2rl+r+r\alpha)/(1+\alpha)}\right]\right)^{(1+\alpha)/4} 
\int_t^T(s-t)^{-\alpha}\left( \E\vert Z^{n,t,x}_s-Z^{k,t,x}_s\vert^2\right)^{(1-\alpha)/4}
\,ds\abs{h}\\
\leq& C\left(1+\vert x \vert\right)^{(2rl+r+r\alpha)/2} 
\left(\int_t^T(s-t)^{-4\alpha/(3+\alpha)}ds\right)^{(3+\alpha)/4}\left( \E\int_t^T\vert Z^{n,t,x}_s-Z^{k,t,x}_s\vert^2\,ds\right)^{(1-\alpha)/4}\abs{h}\\
& \rightarrow 0 \qquad \text{as}\;n,k\rightarrow \infty
\end{align*}
for all $t, x\in\, H$, where in the last passage we have used Proposition \ref{prop-approx-bsde}.
Now we know that for all $t\in\,[0,T[,\,x\in\, H$ the sequences $(v_n(t,x))_n$, and $(\nabla^B v_n(t,x))_n$ converge
and we denote by $\bar v(t,x)$ and $L(t,x)$ respectively their limits.
To conclude we want to show that $\bar v$ is a continuous function, $B$-G\^ateaux differentiable with respect to $x$, $L(t,x)=\nabla^B \bar v (t,x)$,
and $\bar v$ is a mild solution to equation (\ref{Kolmo}). 

\noindent At first we notice that, since
\begin{equation*}
 \vert v^n(t,x) \vert \leq C\left(1+\vert x\vert^{r+1} \right), \qquad
\vert \nabla^B v^n(t,x) \vert \leq C\left(1+\vert x\vert^{r} \right),
\end{equation*}
where $C$ is a constant that does not depend on $n$, $t$, $x$,
then also
\begin{equation*}
 \vert \bar v(t,x) \vert \leq C\left(1+\vert x\vert^{r+1} \right), \qquad
\vert L(t,x) \vert \leq C\left(1+\vert x\vert^{r} \right),
\end{equation*}
where $C$ is the same constant as before.
So, by passing to the limit in (\ref{solmildkolmon}), and also by applying the dominated convergence theorem,
we get 
\begin{equation}
\bar v(t,x)=P_{t,T}\left[  \phi\right]  \left(  x\right)  +\int_{t}^{T}%
P_{t,s}\left[  \psi\left(s,\cdot,\bar v\left(  s,\cdot\right),L\left(  s,\cdot\right)\right)  \right] (  x)  ds,
 \label{solmildkolmolimit}%
\end{equation}
and we can deduce that $\bar v:[0,T]\times H \rightarrow \R$ is a continuous function.
By differentiating (\ref{solmildkolmon}),
we get for all $h \in \Xi$
\begin{align*}
\nabla^B v_n(t,x)h=&\nabla^BP_{t,T}\left[  \phi_n\right]  \left(  x\right) h +\int_{t}^{T}%
\nabla^BP_{t,s}\left[  \psi_n\left(s,\cdot,v_n\left(  s,\cdot\right),\nabla^B 
v_n\left(  s,\cdot\right)\right)  \right] (  x)h  ds\\
 =&\int_{H}\phi_n\left(  z+e^{(T-t)A}x\right)
  \left\langle Q_{T-t}^{-1/2}%
e^{(T-t)A}Bh,Q_{T-t}^{-1/2}y\right\rangle \mathcal{N}\left(  0,Q_{T-t}\right)  \left(
dz\right)\\
&+\int_t^T\int_{H}\psi_n\left(s,  z+e^{(s-t)A}x, v_n(s,  z+e^{(s-t)A}x), \nabla^B v_n(s,  z+e^{(s-t)A}x)\right)\\
  &\times \quad\left\langle Q_{s-t}^{-1/2}%
e^{(s-t)A}Bh,Q_{s-t}^{-1/2}y\right\rangle \mathcal{N}\left(  0,Q_{s-t}\right)  \left(
dz\right)ds.
\end{align*}
By passing to the limit and by applying the dominated convergence theorem,
we get
\begin{align}
\nonumber L(t,x)h=&\int_{H}\phi\left(  z+e^{(T-t)A}x\right)
  \left\langle Q_{T-t}^{-1/2}%
e^{(T-t)A}Bh,Q_{T-t}^{-1/2}y\right\rangle \mathcal{N}\left(  0,Q_{T-t}\right)  \left(
dz\right)\\ 
\nonumber &+\int_t^T\int_{H}\psi\left(s,  z+e^{(s-t)A}x,\bar v(s,  z+e^{(s-t)A}x), L(s,  z+e^{(s-t)A}x)\right)\\
\nonumber &\times\left\langle Q_{s-t}^{-1/2}%
e^{(s-t)A}Bh,Q_{s-t}^{-1/2}y\right\rangle \mathcal{N}\left(  0,Q_{s-t}\right)  \left(
dz\right)ds\\ 
 =&\nabla^BP_{t,T}\left[  \phi\right]  \left(  x\right) h +\int_{t}^{T}%
\nabla^BP_{t,s}\left[  \psi\left(s,\cdot,\bar v\left(  s,\cdot\right),L\left(  s,\cdot\right)\right)  \right] (  x)h  ds. \label{solmildkolmolimit-diffle}
\end{align}
So, in particular we deduce that $L:[0,T)\times H \rightarrow \Xi^{\star}$ is a continuous function.
As a consequence $\psi\left(s,\cdot,\bar v\left(  s,\cdot\right),L\left(  s,\cdot\right)\right) $ is a continuous function,
so by considering (\ref{solmildkolmolimit}) and taking into account the smoothing properties of the transition semigroup $(P_{t,T})_t$,
we deduce that $\bar v:[0,T)\times H \rightarrow \R$ is a $B$-G\^ateaux differentiable function.
Taking the $B$-derivative in (\ref{solmildkolmolimit}) we get for all $h \in \Xi$
\[
\nabla^B\bar v(t,x)h=\nabla^BP_{t,T}\left[  \phi\right]  \left(  x\right)h  +\int_{t}^{T}%
\nabla^BP_{t,s}\left[  \psi\left(s,\cdot,\bar v\left(  s,\cdot\right),L\left(  s,\cdot\right)\right)  \right] (  x)h  ds,
 \]
and by comparing this equation with (\ref{solmildkolmolimit-diffle}) we finally deduce that $\nabla^B \bar v(t,x)=L(t,x)$ and
that $\bar v$ is a mild solution to equation (\ref{Kolmo}). It remains to show that it is the unique mild solution.

\noindent In order to show uniqueness, we  notice that
$\bar v(t,x)=Y_t^{t,x}$, where $Y$ solves the BSDE (\ref{bsde}). 
It remains to show that for every $\tau \in [0,T]$,
$\nabla^{B}\bar v (\tau,X_\tau^{t,x})=Z_\tau^{t,x}$,
where $Z^{t,x}$ is the limit of $Z^{n,t,x}$ in
$L^2(\Omega\times[0,T])$, so in particular $dt \times d\P$-a.s. unless passing to a subsequence.
We already know that for every $n$ $Z_t^{n,t,x}=\nabla^{B}v^n(t,x)$, and
$(\nabla^{B}v_n(t,x))_n$ converges to $\nabla^B \bar v$. Consequently
$\nabla^{B}v^n(\tau,X_\tau^{t,x})\rightarrow \nabla^{B}\bar v(\tau,X_\tau^{t,x})$
$dt \times d\P$-a.s. in $[0,T)\times \Omega$, and
$\nabla^{B}\bar v(\tau,X_\tau^{t,x})=Z_\tau^{t,x}$ $\P$-a.s. for a.a.
$\tau\in[t,T]$.
Since $(Y,Z)$ solves the BSDE (\ref{bsde}), with $Y_t^{t,x}=\bar v(t,x)$,
by previous arguments we get $Z_t^{t,x}=\nabla^{B}\bar v(t,x)$.
By the same arguments of the proof of Theorem \ref{teokolmosmooth},
the solution of the Kolmogorov equation (\ref{Kolmo}) is unique since the solution of the corresponding BSDE
is unique, and this concludes the proof of Theorem \ref{teokolmolip}.
\qed

We now state and prove a theorem analogous to Theorem \ref{teokolmolip}
for the case of a Kolmogorov equation
related to a perturbed Ornstein-Uhlenbeck transition semigroup.

\noindent In the proof of Theorem \ref{teokolmolip} the crucial point is the
regularizing property \ref{ipH su fi} for the Ornstein-Uhlenbeck transition semigroup. We recall that
in \cite{Mas1} regularizing properties of the
Ornstein-Uhlenbeck transition semigroup are linked to regularizing properties of the perturbed
Ornstein-Uhlenbeck transition semigroup related to the process $X^{t,x}$ defined in (\ref{ornstein-pert}). Namely,
in order to verify Hypothesis \ref{ipH su fi}  for the
transition semigroup of the perturbed Ornstein-Uhlenbeck 
process (\ref{ornstein-pert}), we usually assume that $A$
and $B$ satisfy Hypotheses \ref{ornstein inclusione} and \ref{ornstein inclusionestima}.
Then we suppose that $\operatorname{Im} (F)\subset\operatorname{Im}(B)$, namely 
\begin{equation}
 \label{ip-F}
F(t,x)=B G(t,x)
\end{equation}
where $G:[0,T]\times H\rightarrow \Xi$ is bounded and Lipschitz continuous with respect to $x$
uniformly with respect to $t$, and $G\in \calg^{0,1}([0,T]\times H)$.
In such a case it has been proved in \cite{Mas1} that the
perturbed Ornstein-Uhlenbeck process has the same regularizing properties than
the corresponding Ornstein-Uhlenbeck process, i.e. the process defined by (\ref{ornstein-pert}) with $F=0$.

\noindent In the proof of the following theorem we will not use directly this assumption
to get the regularizing property of the perturbed Ornstein-Uhlenbeck transition semigroup,
but an equivalent representation
of the mild solution in terms of an Ornstein-Uhlenbeck transition semigroup. Also in this way, we have to assume that
$F$ satisfies (\ref{ip-F}) as well.

\begin{theorem}\label{teokolmolippert}
Let $A,\, B,\, F$ be the coefficients in the definition of
the perturbed Ornstein-Uhlenbeck process (\ref{ornstein-pert}). Assume that Hypotheses \ref{ip su AB} and \ref{ip-psi-phi}
hold true, and let $F$ satisfy (\ref{ip-F}) with $G\in \calg^{0,1}([0,T]\times H)$ a Lipschitz continuous bounded function.
Moreover assume that the Ornstein-Uhlenbeck transition semigroup defined by (\ref{ornstein-pert}) with $F=0$ 
 satisfies Hypothesis \ref{ipH su fi}. Then, according to Definition \ref{defsolmildkolmo},
equation (\ref{Kolmo}) admits a unique mild solution. 
\end{theorem}
\textbf{Proof.} As already mentioned, in order to prove the theorem for a perturbed Ornstein-Uhlenbeck process, we look
for an equivalent representation of the mild solution in terms of an Ornstein-Uhlenbeck transition semigroup.
To this aim, notice that, at least in the case of $\phi$ and $\psi$ differentiable,
we can apply the Girsanov theorem in the forward-backward system 
\begin{equation*}
\left\{
\begin{array}
[c]{l}%
dX_\tau  =AX_\tau d\tau+BG(\tau,X_\tau) d\tau+BdW_\tau
,\text{ \ \ \ }\tau\in\left[  t,T\right] \\
X_\tau =x,\text{ \ \ \ }\tau\in\left[  0,t\right], \\
\dis
 dY_\tau^{t,x}=-\psi(\tau, X^{t,x}_\tau,Y_\tau^{t,x},Z_\tau^{t,x})\;d\tau+Z^{t,x}_\tau\;dW_\tau,
 \qquad \tau\in [0,T],
  \\\dis
  Y_T^{t,x}=\phi(X_T^{t,x}),
\end{array}
\right.  %
\end{equation*}
or we can follow \cite{Go1}.
We get that the mild solution of equation (\ref{Kolmo}) can be represented, for all $t \in [0,T]$, $x \in H$, as
\begin{equation*}
v(t,x)=R_{t,T}\left[  \phi\right]  \left(  x\right)  +\int_{t}^{T}%
R_{t,s}\left[  \psi\left(s,\cdot,v\left(  s,\cdot\right),\nabla^B 
v\left(  s,\cdot\right)\right)  \right] (  x)  ds +\int_{t}^{T}%
R_{t,s}\left[ \nabla^B v(s\cdot) G (s,\cdot)\right] (  x)  ds. 
\end{equation*}
Here $(R_{t,T})_{t\in[0,T]}$ is the transition semigroup of the corresponding
Ornstein-Uhlenbeck process
\begin{equation*}
\left\{
\begin{array}
[c]{l}%
dX_\tau  =AX_\tau d\tau+BdW_\tau
,\text{ \ \ \ }\tau\in\left[  t,T\right], \\
X_t =x,\text{ \ \ \ }\tau\in\left[  0,t\right].
\end{array}
\right.
\end{equation*}
The new Hamiltonian function is given by 
\begin{equation}\label{newham}
  \tilde \psi(t,x,y,z):=\psi(t,x,y,z)+zG(x)
\end{equation}
and satisfies Hypothesis \ref{ip-psi-phi}. Moreover $G$ by our assumptions is differentiable
so that 
\[
 \tilde \psi_n(t,x,y,z):=\psi_n(t,x,y,z)+zG(x)
\]
where $\psi_n$ is defined in (\ref{infsupconvpsi}).
So we can apply Theorem \ref{teokolmolip}, and the general case of a perturbed Ornstein-Uhlenbeck process is covered.
\qed

\begin{remark}\label{remark-teokolmolippert}
It is possible to show by standard approximations that results stated in Theorem \ref{teokolmolippert} are still true
by taking $G$ only Lipschitz continuous: indeed in this case the new Hamiltonian function
$\tilde\psi$ defined in (\ref{newham}) still satisfies Hypothesis \ref{ip-psi-phi}.
\end{remark}

\section{Application to control}
\label{applic contr}

\subsection{Optimal stochastic control problem}

We formulate the optimal stochastic control problem in the strong
sense. Let $\left(  \Omega,\mathcal{F},\mathbb{P}\right)  $ be a given
complete probability space with a filtration $\left(  \mathcal{F}_{\tau
}\right)  _{\tau\geq0}$ satisfying the usual conditions. $\left\{  W\left(
\tau\right)  ,\tau\geq0\right\}  $ is a cylindrical Wiener process on $H
$\ with respect to $\left(  \mathcal{F}_{\tau}\right)  _{\tau\geq0}$. The
control $u$ is an $\left(  \mathcal{F}_{\tau}\right)  _{\tau}$-predictable
process with values in a closed set $K$ of a normed space $U$; in the following
we will make further assumptions on the control process.
Let us consider the function $R: U\rightarrow H$ and the controlled state equation
\begin{equation}
\left\{
\begin{array}
[c]{l}%
dX^{u}_\tau  =\left[  AX^{u}_\tau+B G(X^u_\tau) +BR\left( u_\tau\right)
 \right]  d\tau+BdW_\tau ,\text{ \ \ \ }\tau\in\left[  t,T\right], \\
X^{u}_t  =x.
\end{array}
\right.  \label{sdecontrolforte}%
\end{equation}
The solution of this equation will be denoted by
$X^{u,t,x}$ or simply by $X^{u}$. $X^u$ is also called
the state, $T>0,$ $t\in\left[  0,T\right]$ are fixed.
The special structure of equation (\ref{sdecontrolforte}) allows to study
the optimal control problem related by means of BSDEs and
(\ref{sdecontrolforte}) leads to a semilinear Hamilton Jacobi Bellman
equation with the structure of the Kolmogorov equation (\ref{Kolmo})
studied in previous sections. The occurrence of the operator
$B$ in the control term is imposed by our techniques, on the contrary the
presence of the operator $R$ allows more generality.

Beside equation (\ref{sdecontrolforte}), we define the cost
\begin{equation}
J\left(  t,x,u\right)  =\mathbb{E}\int_{t}^{T}\left[
\bar g\left(s,X^{u}_s\right) +g\left(u_s\right)\right]ds+\mathbb{E}\phi\left(X^{u}_T\right). 
\label{cost}%
\end{equation}
for real functions $\bar g$ on $[0,T]\times H$, $g$ on $U$ and $\mathbb{\phi}$ on $H$.

\noindent The control problem in strong formulation is to
minimize this functional $J$ over all admissible controls $u$.
We make the following assumptions on the cost $J$.

\begin{hypothesis}
\label{ip costo}

\begin{enumerate}
\item $g:U\rightarrow\R$ is measurable. For some $1<q\leq 2$ there exists a constant $c>0$
such that
\begin{equation}
0\leq g(u)\leq c(1+ \vert u \vert ^q) 
\end{equation}
and there exist $R>0$, $C>0$ such that
\begin{equation}
\label{crescita costo}
g(u)\geq C \vert u \vert ^q \qquad \text{for every }u \in K  \text{ such that } \vert u\vert \geq R.
\end{equation}

\item  There exist $r \in [0,q-1[$, $C>0$, $\alpha>0$ and $\beta>0$ such that  for all $(t,x,x')\in [0,T]\times H\times H$
\[
 \vert \bar g(t,x)-\bar g(t,x')\vert \leq 
\left(C+\frac{\beta}{2}\vert x\vert^r+\frac{\beta}{2}\vert x'\vert^r\right)\vert x -x'\vert;
\]
\[
 \vert \phi(x)-\phi(x')\vert \leq 
\left(C+\frac{\alpha}{2}\vert x\vert^r+\frac{\alpha}{2}\vert x'\vert^r\right)\vert x -x'\vert.
\]

\end{enumerate}

\end{hypothesis}

In the following we denote by 
$\mathcal{A}_{d}$ the set of
admissible controls, that is the $K$-valued predictable processes such that 
\[
 \E \int_0^T \vert u_t \vert ^q dt <+\infty.
\]
This summability requirement is justified  by (\ref{crescita costo}):
a control process which is not $q$-summable would have infinite cost.

\noindent We denote by $J^{\ast}\left(  t,x\right)  =\inf_{u\in\mathcal{A}%
_{d}}J\left(  t,x,u\right)  $ the value function of the problem and, if it
exists, by $u^{\ast}$ the control realizing the infimum, which is called
optimal control.

\noindent We make the following assumptions on $R$.

\begin{hypothesis}
\label{ip aggiuntive}
$R: U\rightarrow H$ is measurable and $\vert R(u)\vert\leq C(1+\vert u\vert) $
for every $u\in U$.
\end{hypothesis}

We have to show that equation (\ref{sdecontrolforte}) admits a unique mild solution, for
every admissible control $u$.
\begin{proposition} \label{prop_sde_contr}
 Let $u$ be an admissible control and assume that Hypothesis \ref{ip su AB} holds true.
Then equation (\ref{sdecontrolforte}) admits a unique mild solution $(X_\tau^u)_{\tau\in[t,T]}$
such that $\E\sup_{\tau\in[t,T]}\vert X_\tau^u\vert^q< \infty$.
\end{proposition}

{\bf Proof.} The proof follows in part the proof of Proposition 2.3 in
 \cite{fuhute}, with some differences since in that paper
 the finite dimensional case is considered and the current cost $g$
 has quadratic growth with respect to $u$, that is to say $q=2$ in (\ref{crescita costo})
 (see also the proof of Proposition 3.16 in \cite{Mas3}, where the case of an Ornstein-Uhlenbeck process is considered). 

 As in \cite{Mas3},  to make an approximation procedure in (\ref{sdecontrolforte})
 we introduce the sequence of stopping times
 $$
 \tau_n=\inf \left\lbrace t\in[0,T]: \E \int_0^t \vert u_s \vert ^q ds>n \right\rbrace \wedge\, T.
 $$
From \cite{fuhute}, we deduce that $\tau_n \rightarrow T$ a.s. in an increasing way as $n \rightarrow +\infty$. Let us define
 \[
  u^n_t=u_t 1_{t\leq\tau_n}+u^0 1_{t>\tau_n},\text{ with }u^0\in K,
 \]
and consider the equation
 \begin{equation}\label{sdecontrolforte-n}
 \left\{
 \begin{array}
 [c]{l}%
 dX^{n}_\tau  =\left[  AX^{n}_\tau+BG\left( X^n_\tau \right) +BR\left( u^n_\tau  \right)
  \right]  d\tau+BdW_\tau ,\text{ \ \ \ }\tau\in\left[  t,T\right], \\
 X^{n}_t  =x.
 \end{array}
 \right.  %
 \end{equation}
 The unique mild solution of equation (\ref{sdecontrolforte-n}) is given by
 \[
  X^n_\tau=e^{(\tau-t)A}x +\int_t^\tau e^{(s-t)A}BG\left( X^n_s \right)
  ds +\int_t^\tau e^{(s-t)A}BR\left( u^n_s \right)
  ds+\int_t^\tau e^{(s-t)A}BdW_s
 \]
 and, by standard calculations, we obtain
 \[
  \E\sup_{\tau\in[t,T]}\vert X^n_\tau\vert^q \leq C\left( \vert x\vert^q+
\E\int_t^T e^{q\omega (s-t)}\vert X^n_s\vert^q\,ds+
 \E\int_t^T (1+\vert u^n_s\vert^q )ds
 +\E\left(\int_t^T e^{2\omega (s-t)}ds\right)^{q/2}\right).
 \]
 Since
 \[
 \E\int_t^T (1+\vert u^n_s\vert^q )ds
 \leq \E\int_t^T (1+\vert u_s\vert^q )ds+T(1+\vert u^0\vert^q)<+\infty,
 \]
 we get, by applying the Gr\"onwall lemma,
that there exists a unique mild solution such that 
\begin{equation}
\label{inegalite uniforme}
\E\left[\sup_{\tau\in[t,T]}\vert X_\tau^n\vert^q \right]\leq C,
\end{equation}
with $C$ that does not depend on $n$.

We have $X_t^n=X_t^{n+1}$ for $t \leq \tau_n$. 
Therefore there exists a process $X$ such that $X_t=X_t^n$ for $t \leq \tau_n$ and $X$ is clearly the required solution. 
The property $\E[\sup_{\tau\in[t,T]}\vert X_\tau\vert^q] <+\infty$ is an immediate consequence of (\ref{inegalite uniforme}).
\qed

We define in a classical way the Hamiltonian function relative to the above
problem:%
\begin{equation*}
h\left(z\right)  =\inf_{u\in K}\left\{  g\left(u\right)
+zR(u)\right\}\quad \forall z\in H .
\end{equation*}
Following the proof of Lemma 3.10 in \cite{Mas3}, we prove that Hypothesis \ref{ip-psi-phi}
is satisfied.
\begin{lemma}
 \label{lemma-hamilton}
Let us define $\psi:[0,T]\times H\times \Xi\rightarrow \R$ by
\begin{equation*}
 \psi(t,x,z):=\bar g(t,x)+h(z)
\end{equation*}
Then $\psi$ satisfies Hypothesis \ref{ip-psi-phi}.
\end{lemma}
{\bf Proof.} The proof follows by our assumptions on $\bar g$ in Hypothesis \ref{ip costo}, and by the proof of Lemma 3.10 in \cite{Mas3}. We notice that the presence of $B G $ in the forward equation can be handled in the
same way as we have done in proposition \ref{prop_sde_contr}, and the polynomial growth of the hamiltonian and of the final condition do not imply substantial changes in the proof.
\qed

\begin{remark}
 We give an example of Hamiltonian we can treat. If in the current cost we take 
$g(u)=\vert u\vert ^q$, $1<q\leq 2$, and in the controlled equation we take $R(u)=u$,
then the Hamiltonian function turns out to be
\[
 \psi(z)=\left( \left(\dfrac{1}{q}\right)^{1/(q-1)}-\left(\dfrac{1}{q}\right)^p\right)\vert z\vert^p 
\]
where $p\geq2$ is the conjugate of $q$.
We underline the fact that our theory covers also the case of
Hamiltonian functions not exactly equal to $\vert z \vert ^p$.
Also notice that the following relation holds true: $l=p-1$,
with $l$ introduced in Hypothesis \ref{ip-psi-phi}.
\end{remark}
We define
\begin{equation}\label{defdigammagrande}
\Gamma(z)=\left\{ u\in U: zR(u)+g(u)=h(z)\right\}.
\end{equation}
If $\Gamma(z) \neq \emptyset$ for every $z\in H$, then by \cite{AuFr} (see Theorems 8.2.10 and
8.2.11), $\Gamma$ admits a measurable selection, i.e. there exists
a measurable function $\gamma: H \rightarrow U$ with
$\gamma(z)\in \Gamma(z)$ for every $z\in H$.

The following theorem
deals with the fundamental relation for the optimal control
by means of backward stochastic differential equations.

\begin{theorem}\label{th-rel-font}
 Assume Hypotheses \ref{ip su AB}, \ref{ipH su fi}, 
\ref{ip costo} and \ref{ip aggiuntive} hold true.
 For every $t\in [0,T]$, $x\in H$ and for all admissible control $u$ we have $J(t,x,u) \geq v(t,x)$, 
 and the equality holds if and only if, for a.a. $s \in [0,T[$, $\mathbb{P}$-a.s.
$$
u_s\in \Gamma\left( \nabla^{B}
v(s ,X^{u,t,x}_s)
\right).
  $$
\end{theorem}
{\bf Proof.} The proof follows the proof of Theorem 3.11 in \cite{Mas3},
with some small mere modifications due to the polynomial growth with respect to $x$ of $v$
and $\nabla^B v$, and due to the presence of $BG$ in the controlled state equation.
 \qed

Under assumptions of Theorem
\ref{th-rel-font}, let us define now the so called
optimal feedback law:
\begin{equation*}
u(s,x)=\gamma\Big(\nabla^{B}
v(s ,X^{u,t,x}_s) \Big),\qquad
s\in [t,T],\;x\in H.
\end{equation*}
Assume that the closed loop equation admits a solution
$\{\overline{X}_s,\;s\in
[t,T]\}$: for all $s \in [0,T]$
\begin{equation}\label{cle}
\overline{X}_s= e^{(s-t)A}x_0
+\int_{t}^s e^{(r-t)A}R(\gamma(\nabla^{B}
v(r ,\overline{X}_r)))\,dr
+\int_t^s e^{(r-t)A}F(\overline{X}_r)\,dr
+\int_t^s e^{(r-t)A}B\,dW_r.
\end{equation}
Then the pair $(\overline{u}=u(s,\overline{X}_s),\overline{X}_s)_{s\in[t,T]}$
is optimal for the control problem.
We notice that existence of a solution of the closed loop
equation is not obvious, due to the lack
of regularity of the feedback law $u$ occurring in
(\ref{cle}).
This problem can be avoided by formulating the optimal control problem in the weak sense, following \cite{FlSo} (see also \cite{fute} and \cite{Mas}).

By an {\em admissible control system}
we mean
$$(\Omega,\mathcal{F},
\left(\mathcal{F}_{t}\right) _{t\geq 0}, \mathbb{P}, W,
u,X^u),$$
where $W$ is an $H$-valued Wiener process, $u$ is an admissible control and $X^u$
solves the controlled equation (\ref{sdecontrolforte}).
The control problem in weak formulation is to minimize
the cost functional over all the admissible control systems.

\begin{theorem}\label{teo su controllo debole}
Assume Hypotheses \ref{ip su AB}, \ref{ipH su fi},
\ref{ip costo} and \ref{ip aggiuntive} hold true.
 For
every $t\in [0,T]$, $x\in H$ and for
 all admissible control systems we have $J(t,x,u)
 \geq v(t,x)$,
  and the
 equality holds if and only if
$$
u_\tau\in \Gamma\left( \nabla^{B}
v(\tau ,X^{u}_\tau)
\right).
  $$
Moreover
assume that the set-valued map $\Gamma$ is non empty and let $\gamma$
be its measurable selection. Then
\begin{equation*}
u_{\tau}=\gamma(\nabla^{B} v(\tau,X^u_\tau ))
,\text{ \ }\mathbb{P}\text{-a.s. for a.a. }\tau
\in\left[ t,T\right],
\end{equation*}
is optimal.

\noindent Finally, the closed loop equation 
\begin{equation*}
 \left\{
\begin{array}
[c]{l}%
dX^{u}_\tau  =\left[  AX^{u}_\tau+BG(X^{u}_\tau) 
+BR\left( \gamma\left(\nabla^{B} v(\tau,X^u_\tau )\right)  \right)
 \right]  d\tau+B dW_\tau ,\text{ \ \ \ }\tau\in\left[  t,T\right], \\
X^{u}_{\tau}  =x ,\text{ \ \ \ }\tau\in\left[  0,t\right].
\end{array}
\right.
\end{equation*}
admits a weak solution
$(\Omega,\mathcal{F},
\left(\mathcal{F}_{t}\right) _{t\geq 0}, \mathbb{P}, W,X)$
which is unique in law and setting
$$
u_{\tau}=\gamma\left(\nabla^{B} v(\tau,X_\tau^u )\right),
$$
we obtain an optimal admissible control system $\left(
W,u,X\right) $.
\end{theorem}

\noindent {\bf Proof.} The proof follows from the fundamental relation
stated in Theorem \ref{th-rel-font}.
The only difference here is the solvability
of the closed loop equation in the weak sense: this is a standard application of the 
Girsanov theorem. Indeed, by Lemma \ref{lemma-hamilton}, see also Lemma 3.10 in \cite{Mas3},
the infimum in the Hamiltonian
is achieved in a ball of radius $C(1+\vert z \vert ^{p-1})$ and so for the optimal control
$u$ the following estimate holds true: $\P$-a.s. and for a.a. $\tau\in[t,T],\;0\leq t\leq T$,
\[
 \vert u_\tau \vert \leq C(1+\vert Z_\tau^{t,x}\vert^{p-1})
=C\left(1+\vert \nabla^{B}v(\tau,X_\tau^{t,x})\vert^{p-1}\right)\leq  
C\left(1+\vert X_\tau^{t,x}\vert^{r(p-1)}\right).
\]
Thanks to this bound and since $r(p-1)<1$, we can apply a Girsanov change of measure
and the conclusion follows in a standard way.
\qed

\begin{remark}\label{rem-fin-contr}
 Notice that in the present section, for the sake of simplicity,
we have considered control problems where the Hamiltonian function
depends only on $\nabla^B v(t,x)$ and not on $v(t,x)$.

\noindent The dependence of the Hamiltonian on the value function is given by taking into account a cost functional of the following form:
\begin{equation*}
J\left(  t,x,u\right)  =\mathbb{E}\int_{t}^{T}\left[\exp\left\lbrace\int_t^s \lambda(u_r)dr\right\rbrace
\bar{g}\left(s,X^{u}_s\right) +g\left(u_s\right)\right]ds+\mathbb{E}\exp\left\lbrace\int_t^T \lambda(u_r)dr\right\rbrace\phi\left(X^{u}_T\right). 
\end{equation*}
In this case the Hamiltonian function is given by
\begin{equation*}
\psi\left(t,x,y,z\right)  =\inf_{u\in K}\left\{  \bar g\left(t,x\right) +g\left(u\right)
+y\lambda(u)+zR(u)\right\}\quad \forall y,z\in H .
\end{equation*}
We also remark that we have focused our attention on a current cost defined by means of $\bar g\left(t,x\right) +g\left(u\right)$,
see (\ref{cost}), in order to verify the assumptions on the Hamiltonian directly thanks to assumptions on $\bar g$ and $g$. We could consider a more general cost given by
\begin{equation*}
J\left(  t,x,u\right)  =\mathbb{E}\int_{t}^{T}
\tilde g\left(s,X^{u}_s,u_s\right) \,ds+\mathbb{E}\phi\left(X^{u}_T\right),
\end{equation*}
and then the Hamiltonian function becomes
\begin{equation*}
\psi\left( t,x,z\right)  =\inf_{u\in K}\left\{ \tilde g\left(t,x,u\right)
+zR(u)\right\}\quad \forall z\in H .
\end{equation*}
Finally we remark that we could also consider a more generic $R$ in equation (\ref{sdecontrolforte}) depending also on $X$
in a Lipschitz continuous way.
\end{remark}

\subsection{Application to a controlled wave equation}
\label{sez-contr-wave}

We can now consider a controlled stochastic wave equation in a complete
probability space $\left(  \Omega,\mathcal{F},\mathbb{P}\right)  $ with a
filtration $\left(  \mathcal{F}_{\tau}\right)  _{\tau\geq0}$ satisfying the
usual conditions. We consider, for $0\leq t\leq\tau\leq T$ and $\xi\in\left[
0,1\right]  $, the following state equation:
\begin{equation}
\left\{
\begin{array}
[c]{l}%
\frac{\partial^{2}}{\partial\tau^{2}}y_\tau\left(\xi\right)  =\frac
{\partial^{2}}{\partial\xi^{2}}y_\tau\left(\xi\right)+f\left(\xi,y_\tau(\xi)\right)  
+u_\tau\left(\xi\right)  +\dot{W}_\tau\left(\xi\right) \\
y_\tau\left(0\right)  =y_\tau\left(1\right)  =0,\\
y_t\left(\xi\right)  =x_{0}\left(  \xi\right)  ,\\
\frac{\partial y_\tau}{\partial\tau}\left( \xi\right)\mid_{\tau=t}  =x_{1}\left(  \xi\right).
\end{array}
\right.  \label{waveequation}%
\end{equation}
$\dot{W}_ \tau\left(\xi\right)  $ is a space-time white noise on $\left[
0,T\right]  \times\left[  0,1\right]  $ and $u_ \tau\left(\cdot\right)  $ is
an admissible control, that is a predictable process
\[
\left(  \Omega,\mathcal{F},\left(  \mathcal{F}_{\tau}\right)  _{\tau\geq
0},\mathbb{P}\right)  \rightarrow L^{2}\left(  0,1\right)  .
\]
Notice that with this square integrability assumption, $u$ satisfies 
the $q$-integrability required in section \ref{applic contr}. 
Moreover we introduce the
cost functional
\[
J\left(  t,x_{0},x_{1},u\right)  =\E\int_{t}^{T}\int_{0}^{1}
\left[\hat{\bar g}\left(  s,\xi,y_s\left(\xi\right)\right) + \hat g\left(u_s\left(\xi\right)  \right)\right]
d\xi ds+\E\int_{0}^{1}\hat\phi\left(  \xi,y_T\left(\xi\right)  \right)
d\xi.
\]
The optimal control problem is to minimize $J$ over all admissible controls.

\begin{hypothesis}
\label{ip costo wave}We make the following assumptions:

\begin{enumerate}

\item $f$ is defined on $\left[  0,1\right]
\times\mathbb{R}  $ and it is measurable. There exists a constant $C>0$ such that, for a.a.  $\xi\in\left[  0,1\right],$
\[
\left|  f\left( \xi,x\right)-f\left(\xi,y\right)  \right|
\leq C\left|  x-y\right|  .
\]
Moreover $f\left(  \xi,\cdot\right)  \in
C^{1}\left(  \mathbb{R}\right)  $.

\item $\hat g:\R\rightarrow\R$ is measurable. For some $1<q\leq 2$ there exist a constant $c>0$
such that
\begin{equation*}
 0\leq \hat g(u)\leq c(1+ \vert u \vert ^q) ,
\end{equation*}
and there exist $R>0,\,C>0$ such that
\begin{equation*}
\hat g(u)\geq C \vert u \vert ^q \qquad \text{for every }u \text{ such that }\vert u\vert \geq R.
\end{equation*}

\item $\hat{\bar g}$ is defined on $\left[  0,T\right]  \times\left[  0,1\right]
\times\mathbb{R}  $ and for a.a. $\tau
\in\left[  0,T\right]  ,$ $\xi\in\left[  0,1\right]  ,$ the map $\hat{\bar g}\left(
\tau,\xi,\cdot\right)  :\R  \rightarrow\R$ is continuous. There exists $r \in [0,q-1[$ such that
for a.a. $\tau\in\left[
0,T\right]  $, $\xi\in\left[ 0,1\right]$ and $x,y \in \mathbb{R}$, 
\[
\left|  \hat{\bar g}\left(  \tau,\xi,x\right)-\hat{\bar g}\left(  \tau,\xi,y\right)  \right|
\leq \left|  x-y\right|\left(C+\dfrac{\beta}{2}\vert x\vert^r+\dfrac{\beta}{2}\vert y\vert^r\right).
\]

\item $\hat\phi:\left[  0,1\right]  \times\mathbb{R}\rightarrow\R$ and for a.a.
$\xi\in\left[  0,1\right]  ,$ $\bar\phi\left(  \xi,\cdot\right)  $ is uniformly
continuous. Moreover there exists $r \in [0,q-1[$ such that for a.a. $\xi\in\left[ 0,1\right]$ and $x,y \in \mathbb{R}$, 
\[
\left| \hat\phi(\xi, x)-\hat\phi(\xi,y)  \right|
\leq \left|  x-y\right|\left(C+\dfrac{\alpha}{2}\vert x\vert^r+\dfrac{\alpha}{2}\vert y\vert^r\right).
\]

\item $x_{0}$, $x_{1}\in L^{2}\left(  \left[  0,1\right]  \right)  $.

\end{enumerate}

\end{hypothesis}

We want to write equation (\ref{waveequation}) in an abstract form. We introduce
the Hilbert space
\[
H=L^{2}\left(  \left[  0,1\right]  \right)  \oplus\mathcal{D}\left(
\Lambda^{-\frac{1}{2}}\right)  =L^{2}\left(  \left[  0,1\right]  \right)
\oplus H^{-1}\left(  \left[  0,1\right]  \right)  .
\]
In fact in the stochastic case the controlled wave equation does not evolve in
$H_{0}^{1}\left(  \left[  0,1\right]  \right)  \oplus L^{2}\left(  \left[
0,1\right]  \right)  $, see \cite{DP1} and also \cite{Mas}. On $H$ we define the operator $A$ by
\[
\mathcal{D}\left(  A\right)  =H_{0}^{1}\left(  \left[  0,1\right]  \right)
\oplus L^{2}\left(  \left[  0,1\right]  \right)  ,\text{ \ \ \ \ }A\left(
\begin{array}
[c]{c}%
y\\
z
\end{array}
\right)  =\left(
\begin{array}
[c]{cc}%
0 & I\\
-\Lambda & 0
\end{array}
\right)  \left(
\begin{array}
[c]{c}%
y\\
z
\end{array}
\right)  ,\text{ \ for every }\left(
\begin{array}
[c]{c}%
y\\
z
\end{array}
\right)  \in\mathcal{D}\left(  A\right)  .
\]
We also set $G:H\rightarrow L^2([0,1])$
\[
 G\left(\left(\begin{array}
[c]{c}%
y\\
z
\end{array}\right) \right)\left(\xi\right):=f(\xi, y(\xi))
\]
for all $\left(\begin{array}
[c]{c}%
y\\
z
\end{array} \right)\in H$ and 
 $B:L^{2}\left(  \left[  0,1\right]  \right)  \longrightarrow H$
with $Bu=\left(
\begin{array}
[c]{c}%
0\\
u
\end{array}
\right)  =\left(
\begin{array}
[c]{c}%
0\\
I
\end{array}
\right)  u$, $u\in L^2([0,1])$.  Thanks to Hypothesis \ref{ip costo wave}, point 1, $F:=BG$ satisfies
Hypothesis \ref{ip su AB}, point 2.

\noindent Equation (\ref{waveequation}) can be rewritten in an abstract way as an
equation in $H$ of the following form:%
\begin{equation}
\left\{
\begin{array}
[c]{l}%
dX^{u}_\tau =AX^{u}_\tau  d\tau+BG\left( X^{u}_\tau\right)d\tau+Bu_\tau d\tau
+BdW_\tau  ,\text{ \ \ \ }\tau\in\left[
t,T\right] \\
X^{u}_t =x,
\end{array}
\right.  \label{waveeqabstract}%
\end{equation}
We notice that by \cite{Mas}, section 6.1, the transition semigroup of the linear uncontrolled wave equation, i.e equation
\ref{waveeqabstract} with $F=0$ and without control, satisfies Hypothesis \ref{ipH su fi}
with $\alpha=1/2$, and by  \cite{Mas1}, the transition semigroup of the uncontrolled wave equation, i.e equation
\ref{waveeqabstract} without control, also satisfies Hypothesis \ref{ipH su fi}
with $\alpha=1/2$.

Moreover, for all $x=\left(\begin{array}
[c]{c}%
y\\
z
\end{array}\right)  \in H$ and for all $u\in L^2([0,1])$, we set
\[%
\begin{array}
[c]{l}%
\bar g\left(  \tau,x\right)  =\left(
{\displaystyle\int_{0}^{1}}
\hat{\bar g}\left(  \tau,\xi,y\left(  \xi\right)   \right)
d\xi\right) ,\qquad g(u) =\left({\displaystyle\int_{0}^{1}}
\hat g\left( u(\xi)
\right)d\xi \right),\\
\phi\left(  x \right) =\left(
{\displaystyle\int_{0}^{1}}
\hat\phi\left(  \xi,y\left(  \xi\right)   \right)
d\xi\right)  .
\end{array}
\] 
Due to the fact that $r<1$ and $q<2$,
it is standard to show that $\bar g$, $g$ and $\phi$ satisfy Hypothesis
\ref{ip costo}.

\noindent In abstract formulation, the cost functional can be written as%
\[
J\left(  t,x,u\right)  =\mathbb{E}\int_{t}^{T}\left(\bar g\left(  s,X^{u}_s\right)+g(u_s)\right)  ds
+\E\phi\left(  X^{u}_T\right)  .
\]
We solve the control problem in its weak formulation, which allows
to make the synthesis of the optimal control by solving the closed loop equation in weak sense.
We define $v$ as the solution of the Hamilton Jacobi Bellman equation associated
to the uncontrolled wave equation.
\begin{theorem}\label{teo su controllo debole wave}
Assume Hypothesis \ref{ip costo wave} holds true.
 For
every $t\in [0,T]$, $x\in H$ and for
 all admissible control systems we have $J(t,x,u(\cdot))
 \geq v(t,x)$,
  and the
 equality holds if and only if
$$
u_\tau\in \Gamma\left( \nabla^{B}
v(\tau ,X^{u}_\tau)
\right),
  $$
where $\Gamma$ has been defined in (\ref{defdigammagrande}).
Moreover
assume that the set-valued map $\Gamma$ is non empty and let $\gamma$
be its measurable selection, then
\begin{equation*}
u_{\tau}=\gamma(\nabla^{B} v(\tau,X^u_\tau ))
,\text{ \ }\mathbb{P}\text{-a.s. for a.a. }\tau
\in\left[ t,T\right]
\end{equation*}
is optimal.

\noindent Finally, the closed loop equation 
\begin{equation*}
 \left\{
\begin{array}
[c]{l}%
dX^{u}_\tau  =\left[  AX^{u}_\tau+BG(X^{u}_\tau) 
+B \gamma\left(\nabla^{B} v(\tau,X^u_\tau )\right) 
 \right]  d\tau+B dW_\tau ,\text{ \ \ \ }\tau\in\left[  t,T\right] \\
X^{u}_t  =x.
\end{array}
\right.
\end{equation*}
admits a weak solution
$(\Omega,\mathcal{F},
\left(\mathcal{F}_{t}\right) _{t\geq 0}, \mathbb{P}, W,X)$
which is unique in law and setting
$$
u_{\tau}=\gamma\left(\nabla^{B} v(\tau,X_\tau^u )\right),
$$
we obtain an optimal admissible control system $\left(
W,u,X\right) $.
\end{theorem}

\noindent {\bf Proof.} The proof follows from Theorem \ref{teo su controllo debole}
by noticing that Hypothesis \ref{ip su AB} and  \ref{ip costo}
follow by Hypothesis \ref{ip costo wave}, Hypothesis \ref{ip aggiuntive} is satisfied since $R$
equals the identity, and as previously noticed 
Hypothesis \ref{ipH su fi} is satisfied by the transition semigroup of the uncontrolled wave equation with $\alpha=1/2$.
\qed

\end{document}